\documentclass{article}
\usepackage{mystyle}
\begin{document}
\title{Measurable cardinals and choiceless axioms}
\author{Gabriel Goldberg\footnote{
    The author was partially supported by NSF Award DMS-1902884. 
}}
\maketitle
\begin{abstract}
    Kunen refuted the existence of an elementary embedding from the universe
    of sets to itself assuming the Axiom of Choice.
    This paper concerns the ramifications of this hypothesis when the Axiom of Choice
    is not assumed. For example, the existence of such an 
    embedding implies that there is a proper class of
    cardinals \(\lambda\) such that \(\lambda^+\) is measurable.
\end{abstract}
\section{Introduction}
\subsection{The Kunen inconsistency}
One of the most influential ideas in the
history of large cardinals is Scott's reformulation
of measurability in terms of elementary embeddings
\cite{Scott}:
the existence of a measurable cardinal
is equivalent to the existence of 
a nontrivial elementary embedding from the universe of sets \(V\)
into a transitive submodel \(M\).
In the late 1960s, Solovay and Reinhardt realized that by imposing
stronger and stronger closure constraints on the model \(M\),
one obtains stronger and stronger large cardinal axioms,
an insight which rapidly led to the discovery of most of the modern large cardinal
hierarchy. Around this time, Reinhardt formulated 
the ultimate large cardinal principle of this kind:
there is an elementary embedding from the universe of sets to itself.\footnote{
    Of course, the identity is such an elementary embedding. Whenever
    we write ``elementary embedding," we will really mean ``nontrivial
    elementary embedding.''
}
Soon after, however, Kunen \cite{Kunen} showed that
this principle is inconsistent:
\begin{thm*}[Kunen]
    There is no elementary embedding from the universe of sets to itself.
\end{thm*}

Kunen's proof relies
heavily on the Axiom of Choice, however, and the question
of whether this is necessary immediate arose.\footnote{The question was first raised
by the anonymous referee of Kunen's paper.}
Decades later, Woodin returned to this question and 
discovered that although the traditional large cardinal hierarchy
stops short at Kunen's bound,
there lies beyond it a further realm of
large cardinal axioms incompatible with the Axiom of Choice,
axioms so absurdly strong
that Reinhardt's so-called ultimate axiom appears tame by comparison.
Yet since their discovery, despite significant efforts of many researchers,
no one has managed to prove the inconsistency of a single 
one of these choiceless large cardinal axioms.
``The difficulty,'' according to Woodin \cite{WoodinTransfinite},
``is that without the Axiom of Choice it is extraordinarily difficult to prove anything about sets.''

One remedy to this difficulty, proposed by Woodin himself \cite[Theorem 227]{SEM},
is to \textit{simulate} the Axiom of Choice using auxiliary large cardinal hypotheses,
especially extendible cardinals. 
Cutolo \cite{Cutolo} expanded on this idea to establish the
striking result that the successor of a singular Berkeley limit of extendible cardinals is measurable.
While the strength of Cutolo's large cardinal hypothesis far surpasses that of
Reinhardt's ultimate axiom,
Asper\'o \cite{AsperoNote} showed that a Reinhardt cardinal alone
implies the existence of elementary embeddings reminiscent
of those associated with extendible cardinals.
Combining Woodin, Cutolo, and Asper\'o's ideas, we show here
that one can simulate the Axiom of Choice using
large cardinal notions 
that follow from Reinhardt's principle.
This allows us to establish some vast generalizations of Cutolo's results.
For example, from Reinhardt's axiom alone, one can establish the
existence of a proper class of measurable successor cardinals.
We are optimistic that these ideas will bring clarity 
to the question of the consistency
of the choiceless hierarchy.

\subsection{Main results}
Gitik \cite{GitikSingular}
showed that it is consistent with ZF that there are no regular uncountable cardinals.
On the other hand, 
if there is an elementary embedding \(j\) from the universe of sets
to itself, its critical point \(\kappa\) is measurable and hence regular.\footnote{
    The critical point of \(j : V\to M\), denoted \(\crit(j)\), is the least ordinal \(\alpha\)
    such that \(j(\alpha) > \alpha\). Scott \cite{Scott} showed that \(\crit(j)\) is measurable.
}
By elementarity, the cardinals \(\kappa_1(j) = j(\kappa)\), \(
    \kappa_2(j) = j(j(\kappa))\), \(\kappa_3(j) = j(j(j(\kappa)))\), and so on
are all regular as well.
Asper\'o asked whether there must be any regular cardinals larger than
their supremum \(\kappa_\omega(j) = \sup_{n < \omega} \kappa_n(j)\).
The first theorem of this paper answers his question positively:
\begin{repthm}{thm:regulars}
Suppose there is an elementary embedding from the universe of sets to itself.
Then there is a proper class of regular cardinals.
\end{repthm}
This theorem is a consequence of the \textit{wellordered collection lemma},
a weak choice principle
that we show follows from choiceless cardinals:\footnote{A Scott cardinal 
is essentially an equivalence class
of the equinumerosity relation; see \cref{section:cardinals}.
}
\begin{repthm}{thm:wellordered_collection}
    Suppose there is an elementary embedding from the universe of sets to itself.
    For every ordinal \(\kappa\), there is a 
    Scott cardinal \(\mathfrak{a}\) such that for any sequence 
    \(\langle A_\alpha : \alpha < \kappa\rangle\) of nonempty
    sets, there is a set \(\sigma\) of cardinality at most \(\mathfrak{a}\) such that 
    \(A_\alpha\cap \sigma\neq \emptyset\) for all \(\alpha < \kappa\).
\end{repthm}

As a consequence of the wellordered collection lemma, we
obtain one of the key combinatorial consequences of the Axiom of Choice:
\begin{repthm}{thm:completeness}
    Suppose there is an elementary embedding from the universe of sets to itself.
    Then for all cardinals \(\kappa\), for all sufficiently large regular cardinals 
    \(\delta\), the closed unbounded filter on \(\delta\) is \(\kappa\)-complete.
\end{repthm}

Having answered Asper\'o's question, it is natural to wonder whether
Reinhardt's principle in fact implies the existence a proper class of \textit{measurable} 
cardinals.
Given Cutolo's result, one would also like to know whether any of these measurable
cardinals are successor cardinals.
\begin{repthm}{thm:measurables}
    Suppose there is an elementary embedding from the universe of sets to itself.
    Then for a closed unbounded class of cardinals \(\kappa\), either \(\kappa\) or \(\kappa^+\)
    is measurable.
\end{repthm}
In particular, for every regular cardinal \(\gamma\), there are arbitrarily
large cardinals \(\kappa\) of cofinality \(\gamma\) such that \(\kappa^+\) is measurable.

There is really only one other principle that is known to imply the existence
of measurable successor cardinals: the Axiom of Determinacy (AD). Solovay
showed
that under AD, \(\aleph_1\) is measurable. Moreover,
there is a unique normal ultrafilter on \(\aleph_1\):
the closed unbounded filter.
Later, Martin showed that \(\aleph_2\) is measurable,
and finally Kunen \cite{KechrisProjective} showed:
\begin{thm*}[Kunen]
Assuming the Axiom of Determinacy,
\(\boldsymbol{\delta}^1_n\) is measurable for all \(n < \omega\).
\end{thm*}
Here \(\boldsymbol{\delta}^1_n\) is the \(n\)-th
\textit{projective ordinal,} the supremum of
all \(\boldsymbol{\Delta}^1_n\)-definable prewellorders of the real
numbers; for example,
\(\boldsymbol{\delta}^1_1 = \aleph_1\) and under AD,
\(\boldsymbol{\delta}^1_2 = \aleph_2\) and
\(\boldsymbol{\delta}^1_3 = \aleph_{\omega+1}\).
Results
of Kechris, Kunen, and Martin imply that all the projective ordinals are successor
cardinals and
that the \(\omega\)-closed unbounded filter on \(\boldsymbol{\delta}^1_n\)
is a normal ultrafilter, generalizing Solovay's theorem.

The next part of the paper attempts
to generalize some of the determinacy 
theory of ultrafilters and closed unbounded sets
to the context of choiceless large cardinal axioms.
If \(F\) is a filter, a set is \textit{\(F\)-positive}
if it intersects every set in \(F\), and an \(F\)-positive set
\(S\) is an \textit{atom of \(F\)}
if \(F\cup \{S\}\) generates an ultrafilter. A filter \(F\) is \textit{atomic}
if every \(F\)-positive set contains an atom.
Assuming AD,
Kechris-Kleinberg-Moschovakis-Woodin
\cite{KechrisKleinbergWoodin} proved the existence of many 
cardinals with the strong partition property
and showed that the closed unbounded filter on such a cardinal is 
atomic.\footnote{
    A cardinal \(\delta\) has the \textit{strong partition property} if 
    any subset of \([\delta]^\delta\)
    either contains or is disjoint from 
    a set of the form \([H]^\delta\)
    for some \(H\in [\delta]^\delta\).
} 

In this situation, one can actually
classify the atoms of the closed unbounded filter in terms of stationary reflection.
A stationary set \(S\) \textit{reflects stationarily in \(T\)} if there is a stationary set
of \(\alpha\in T\) such that \(S\cap \alpha\) is stationary;
\(S\) \textit{reflects fully in \(T\)} if for all but a nonstationary
set of \(\alpha\in T\), \(S\cap \alpha\) is stationary.
A stationary set is \textit{thin} if it does not reflect stationarily in itself.
Since every stationary set contains a thin stationary set,
every atom of the closed unbounded filter is thin. Assuming the strong partition property,
Kechris-Kleinberg-Moschovakis-Woodin proved the converse:
every thin stationary subset of
a strong partition cardinal is an atom of the closed unbounded filter.

It is not possible to prove such a strong theorem from choiceless
large cardinal axioms: for example,
an easy forcing argument shows that it is consistent
with choiceless axioms that for all regular cardinals \(\delta\), 
the \(\omega\)-closed unbounded filter is \textit{not} an ultrafilter,
which means that the thin set 
\(S^\delta_\omega = \{\alpha < \delta : \cf(\alpha) = \omega\}\)
is not an atom of the closed unbounded filter.\footnote{
    Just let \(G\subseteq \text{Col}(\omega,\omega_1)\) be \(V\)-generic.
    Then for all ordinals \(\delta\), 
    \((S^\delta_\omega)^{V[G]}\) is the union of
    \((S^\delta_\omega)^V\) and \((S^\delta_{\omega_1})^V\), each of which is stationary
    in \(V[G]\) assuming \(\cf^V(\delta)\geq \omega_2\).
}
This motivates the following theorem, which says
that thin sets are \textit{almost} atoms:
\begin{repthm}{thm:atoms}
    Suppose \(j : V\to V\) is an elementary embedding.
    Then for all sufficiently large regular cardinals \(\delta\),
    any thin stationary subset of \(\delta\)
    can be written as the disjoint union of at most \(\kappa_\omega(j)\)
    atoms of the closed unbounded filter.
\end{repthm}

Jech's \textit{reflection order} on stationary sets
is defined by setting
\(S < T\) if \(S\) reflects fully in \(T\). 
Steel proved that if \(\delta\) is a strong partition
cardinal, then the reflection order wellorders
the atoms of the closed unbounded filter on \(\delta\)
modulo the nonstationary ideal.
Assuming Reinhardt's principle, one
can prove an analog of Steel's result for all sufficiently large regular
cardinals:
\begin{repthm}{thm:reflection_order}
    If \(j : V\to V\) is an elementary embedding, 
    then for all sufficiently large regular cardinals
    \(\delta\), any set of reflection incomparable
    atoms of the closed unbounded filter
    has cardinality less than or equal to \(\kappa_\omega(j)\).
\end{repthm}

The proofs of these theorems involve an order on ultrafilters
known as the \textit{Ketonen order}, whose tortuous history we
now describe. The discovery of the Ketonen order
was precipitated by Kunen's construction 
of a normal ultrafilter that concentrates
on nonmeasurable cardinals. 
Ketonen \cite{Ketonen}, then completing his dissertation under Kunen,
realized that implicit in
this proof was a natural order on weakly normal ultrafilters.
He and Kunen collaborated to prove the wellfoundedness of this order,
which Ketonen needed to prove his celebrated result
that if every regular cardinal above \(\kappa\) carries
a \(\kappa\)-complete uniform ultrafilter, then \(\kappa\) is strongly
compact (answering a question of Kunen).

The next year, also inspired by Kunen's construction, Mitchell \cite{Mitchell} independently 
discovered Ketonen's order, or rather its
restriction to normal ultrafilters.
Mitchell proved that this order is linear in canonical inner models
of large cardinal axioms. Since then, the \textit{Mitchell order} has 
become a fundamental object of study in large cardinal theory. 

Ketonen's order, on the other hand, seems to have been 
forgotten completely until almost half a century later, 
the author independently discovered a generalization of his
order to \textit{all} countably complete ultrafilters on ordinals \cite{KO}:
if \(U\) and \(W\) are countably complete
ultrafilters on ordinals, set \(U\ke W\) if
there is a sequence ultrafilters \(U_\alpha\) on \(\alpha\), 
defined for \(W\)-almost all ordinals \(\alpha\),
such that \(A\in U\) if and only if \(A\cap \alpha\in U_\alpha\) for
\(W\)-almost all \(\alpha\). Probably Ketonen thought about this
generalization, but the Ketonen-Kunen wellfoundedness proof does not generalize to
arbitrary ultrafilters. It turns out, however, that
\(\ke\) \textit{is} wellfounded, 
but this requires an argument that is completely different
from Ketonen and Kunen's.
Moreover, like the Mitchell order, \(\ke\) is linear in all known canonical inner models
of large cardinal axioms, although again the proof of this 
is completely different from Mitchell's.

The order \(\ke\) is now known as the \textit{Ketonen order.}
In the context of the Axiom of Choice, the linearity of the Ketonen order is equivalent to 
the Ultrapower Axiom \cite{UA}, a principle with many consequences
in large cardinal theory,
but in this paper, we will apply linearity
properties of the Ketonen order to the theory of choiceless cardinals.
The key phenomenon is that choiceless cardinals imply that the Ketonen order is \textit{almost} linear:
\begin{repthm}{thm:antichains}
    Suppose \(j : V\to V\) is an elementary embedding.
    Then for all sufficiently large cardinals \(\kappa\),
    any set of Ketonen incomparable \(\kappa\)-complete
    ultrafilters has cardinality less than or equal to \(\kappa_\omega(j)\).
\end{repthm}
In other words, choiceless large cardinal axioms almost imply the Ultrapower Axiom. 
One cannot hope to prove such a theorem
from traditional large cardinal
axioms. It is open whether something analogous follows from the Axiom of Determinacy.

The semi-linearity of the Ketonen order is also a component of the proof of one
of the main theorems of this paper. 
In the context of ZFC, 
an uncountable cardinal 
\(\kappa\) is said to be \textit{strongly compact} if every \(\kappa\)-complete filter
extends to a \(\kappa\)-complete ultrafilter. 
Kunen realized, however, that 
this \textit{filter extension property} occurs in nature:
\begin{thm*}[Kunen]
    Assume \(\AD + \textnormal{DC}_{\mathbb R}\).
    If \(\beta\) is an ordinal that
    is the surjective image of \(\mathbb R\),
    then every countably complete filter on \(\beta\)
    extends to a countably complete ultrafilter.
\end{thm*}
In order to prove a version of Kunen's theorem 
from choiceless large cardinal axioms,
we extend the Ketonen order to a wellfounded partial order on filters. 
We then prove the filter extension property
by induction on this order.
\begin{repthm}{thm:filter_extension}
    Assume there is an elementary embedding from the universe of sets to itself.
    Then for a closed unbounded class of cardinals \(\kappa\),
    every \(\kappa\)-complete filter on an ordinal extends to a \(\kappa\)-complete ultrafilter.
\end{repthm}
\section{Choicelike consequences of choiceless axioms}
\subsection{Cardinality without choice}\label{section:cardinals}
We write \(|X| \leq |Y|\) to mean that there is an injection from \(X\)
to \(Y\), \(|X|\leq^*|Y|\) to mean that there is a partial surjection
from \(Y\) to \(X\),
and \(|X| = |Y|\) to mean there is a bijection between \(X\) and \(Y\). 
The \textit{Scott rank of \(X\)}, 
denoted by \(\scott(X)\), is the least ordinal \(\alpha\) such that
\(|X|\leq |V_\alpha|\);
the \textit{dual Scott rank of \(X\)}, 
denoted by \(\dscott(X)\), is the least ordinal \(\alpha\)
such that \(|X| \leq^* |V_\alpha|\).
Notice that \(\dscott(X) \leq \scott(X) \leq \dscott(X)+1\).

The \textit{Scott cardinality of \(X\)}, denoted by \(|X|\), 
is the family of sets \(Y\) of rank \(\scott(X)\)
such that \(|X| = |Y|\).
A \textit{Scott cardinal} is a set that is the Scott cardinality of some set, 
while a \textit{cardinal}, as usual, is an ordinal that is not in bijection with
any of its predecessors. If \(X\) is wellorderable, we abuse notation by letting
\(|X|\) denote the unique cardinal in bijection with \(X\)
rather than the Scott cardinality of \(X\).

\subsection{Almost extendibility and supercompactness}
A cardinal \(\lambda\) is \textit{rank Berkeley}
if for all \(\alpha < \lambda \leq \beta\),
there is an elementary embedding \(j : V_\beta\to V_\beta\)
with \(\alpha < \crit(j) < \lambda\).
First introduced by Schlutzenberg, 
rank Berkeley cardinals are a weakening of 
Reinhardt cardinals
that have the advantage of being first-order definable.
With the exception of the following proposition, we will 
have nothing to say about proper classes and Reinhardt cardinals,
and instead, we will consider rank Berkeley cardinals in 
first-order ZF.
\begin{prp}[NBG]
    If \(j : V\to V\) is a nontrivial elementary embedding,
    then \(\kappa_\omega(j)\) is rank Berkeley.
    \begin{proof}
        Let \(\lambda = \kappa_\omega(j)\).
        Assume the proposition fails, and consider the least \(\alpha\)
        such that there is no elementary \(k : V_\alpha\to V_\alpha\)
        such that \(\kappa_\omega(k) = \lambda\).
        Then \(\alpha\) is definable from \(\lambda\),
        and hence \(j(\alpha) = \alpha\).
        But the embedding \(k = j\restriction V_\alpha\) contradicts 
        the definition of \(\alpha\). 
    \end{proof}
\end{prp}
The rest of this section deduces the existence
of many pseudo-extendible and pseudo-supercompact cardinals
from the existence of a rank Berkeley cardinal.
This is similar to work of Asper\'o \cite{AsperoNote}.

For any set \(X\), let \(\theta(X)\) be the least ordinal
that is not the range of a function on \(X\).
\begin{prp}\label{prp:bounded_vopenka}
    Suppose \(\alpha\) is a limit ordinal and
    \(j : V_{\alpha}\to V_{\alpha}\) is elementary.
    Suppose \(\mathcal S\) is a subset of \(V_\alpha\) 
    consisting of structures in a fixed finite signature,
    \(j(\mathcal S) = \mathcal S\),\footnote{More formally, 
    we mean that \(\bigcup_{x\in V_\alpha} j(\mathcal S\cap x) = \mathcal S\).} and 
    \(\theta(\mathcal S) \geq \kappa_\omega(j)\). Then
    there exist distinct \(M_0\) and \(M_1\) in \(\mathcal S\)
    such that \(M_0\) elementarily embeds into \(M_1\).
    \begin{proof}
        Suppose not. Let \(\kappa = \crit(j)\) and let \(\lambda=\kappa_\omega(j)\). 
        By replacing \(j\) with \(j(j)\) if necessary, we can assume without loss of generality
        that \(\cf(\alpha) \neq \kappa\). (This was pointed out by the anonymous referee.)

        For \(\beta < \alpha\), let \(\mathcal S_\beta = \mathcal S\cap V_\beta\).
        We claim that for some \(\beta < \alpha\), 
        \(\theta(\mathcal S_\beta) > \kappa\).
        Let \(f : \mathcal S\to \kappa\) be a surjection.
        For each \(\beta < \alpha\), let \(\xi_\beta = \ot (f[\mathcal S_\beta])\).
        Then \(\sup_{\beta < \alpha} \xi_\beta = \kappa\). 
        Since \(\cf(\alpha) \neq \kappa\),
        \(\xi_\beta = \kappa\) for some \(\beta < \alpha\). It follows that
        \(\theta(\mathcal S_\beta) > \kappa\).

        Fix \(\beta\) such that \(\theta(\mathcal S_\beta) > \kappa\).
        Let \(g : \mathcal S_\beta\to \kappa\) be a surjection.
        Note that \(j(g)\) is a surjection from \(j(\mathcal S_\beta)\) to
        \(j(\kappa)\), while \(j(g)(j(M)) = j(g(M)) = g(M) < \kappa\) for all
        \(M\in \mathcal S_\beta\). It follows that there is some \(M_0\in j(\mathcal S_\beta)\)
        such that \(M_0\notin j[\mathcal S_\beta]\). As a consequence, letting \(M_1 = j(M_0)\),
        we have \(M_0\neq M_1\). Since \(j(\mathcal S) = \mathcal S\),
        \(M_1\in \mathcal S\). Moreover, \(j\) restricts to an elementary embedding
        from \(M_0\) to \(M_1\).
    \end{proof}
\end{prp}
The \textit{lightface Vop\v{e}nka principle} states that for all 
parameter-free definable classes \(\mathcal S\) of structures in a fixed finite
signature, there exist distinct structures \(M_0\) and \(M_1\) in \(\mathcal S\)
such that \(M_0\) elementarily embeds into \(M_1\). 
The lightface
Vop\v{e}nka principle implies Vop\v{e}nka's principle for any class of 
structures definable using ordinal parameters. (Consider the least counterexample.)
\begin{cor}
    If \(\lambda\) is rank Berkeley, then for any ordinal definable set of structures
    \(\mathcal S\) in a finite signature
    such that \(\theta(\mathcal S)\geq \lambda\),
    there exist distinct structures \(M_0\) and \(M_1\) in \(\mathcal S\)
    such that \(M_0\) elementarily embeds into \(M_1\).
    As a consequence, the lightface Vop\v{e}nka principle holds.\qed
\end{cor}
A cardinal \(\eta\) is \textit{\((\gamma,\infty)\)-extendible}
if for all \(\nu > \eta\),
there is an elementary embedding \(\pi : V_\nu\to V_{\nu'}\)
such that \(\pi(\eta) > \nu\) and \(\pi(\gamma) = \gamma\). 
\begin{lma}\label{lma:lightface_to_extendible}
    Assume the lightface Vop\v{e}nka principle.
    Then for all ordinals \(\gamma\), there is a \((\gamma,\infty)\)-extendible cardinal.
    \begin{proof}
        Assume not.
        Define a continuous sequence of ordinals \(\langle\eta_\xi :\xi \in \Ord\rangle\)
        by transfinite recursion,
        letting \(\eta_{\xi+1}\) be the least ordinal \(\nu > \eta_\xi\) such that
        there is no \(\pi : V_{\nu}\to V_{\nu'}\) such that \(\pi(\eta_\xi) \geq \nu\).
        Let \(\mathcal S\) be the class of structures
        \(\mathcal M_\xi = (V_{\eta_{\xi+1}},\eta_\xi,\gamma)\). 
        Applying the lightface Vop\v{e}nka principle to the ordinal definable
        class \(\mathcal S\), we obtain ordinals \(\xi_0 < \xi_1\)
        and an elementary embedding \(\pi : \mathcal M_{\xi_0}\to \mathcal M_{\xi_1}\). 
        This means \(\pi : V_{\eta_{\xi_0+1}}\to V_{\eta_{\xi_1+1}}\) 
        is elementary, \(\pi(\eta_{\xi_0}) = \eta_{\xi_1}\),
        and \(\pi(\gamma) = \gamma\).
        This contradicts the definition of 
        \(\eta_{\xi_0+1}\) since 
        \(\pi(\eta_{\xi_0}) = \eta_{\xi_1} \geq \eta_{\xi_0+1}\).
    \end{proof}
\end{lma}

A cardinal \(\eta\) is \textit{almost extendible}
if it is \((\gamma,\infty)\)-extendible for all \(\gamma <\eta\).
\begin{cor}\label{cor:closed unbounded_extendible}
    Assume the lightface Vop\v{e}nka principle.
    Then there is a closed unbounded class of
    almost extendible cardinals.\qed
\end{cor}

A cardinal \(\eta\) is \textit{\((\gamma,\infty)\)-supercompact}
if for all \(\nu > \eta\),
for some \(\bar \nu < \eta\),
there is an elementary embedding \(\pi : V_{\bar \nu}\to V_\nu\)
such that \(\pi(\gamma) = \gamma\).
Note that if \(\eta\) is \((\gamma,\infty)\)-supercompact then 
for all \(\nu > \eta\) and all \(x\in V_\nu\),
for some \(\bar \nu < \eta\),
there is an elementary embedding \(\pi : V_{\bar \nu}\to V_\nu\)
such that \(\pi(\gamma) = \gamma\) and \(x\in \ran(\pi)\).\footnote{To see this, 
    first consider the case that \(x\) is an ordinal \(\alpha < \nu\), which is
    easy to handle since one can find \(\pi : V_{\bar \nu'}\to V_{\nu'}\) where \(\nu'\) is an ordinal coding the pair
    \(\langle \nu ,\alpha\rangle\). To handle a general \(x\), take \(\bar \nu < \eta\) and 
    \(\pi : V_{\bar \nu+1}\to V_{\nu+1}\)
    with \(\eta\in \ran(\pi)\). Then note that for all \(x\in \pi[V_{\bar \nu+1}]\),
    if \(x\in V_\nu\), then for some
    \(\bar \nu' < \eta\), there is an elementary
    embedding \(\pi' : V_{\bar \nu'}\to V_{\nu}\)
    such that \(x\in \ran(\pi')\); namely, \(\pi'=\pi\restriction V_{\bar \nu}\). 
    Since \(\pi[V_{\bar \nu+1}]\preceq V_{\nu + 1}\),
    it follows that for all \(x\in V_{\nu}\), 
    for some
    \(\bar \nu' < \eta\), there is an elementary
    embedding \(\pi' : V_{\bar \nu'}\to V_{\nu}\)
    such that \(x\in \ran(\pi')\), as desired. \label{foot:x_in_ran}
}
A cardinal
\(\eta\) is \textit{almost supercompact} if it is
\((\gamma,\infty)\)-supercompact for all \(\gamma < \eta\).
\begin{lma}\label{lma:extendible_to_sc}
    \begin{enumerate}[(1)]
        \item Every almost extendible cardinal is almost supercompact.\label{item:ex_to_sc}
        \item If \(\eta\) is almost supercompact, \(\varphi\) is a formula, and \(\gamma < \eta\),
        then the least ordinal \(\alpha\) such that \(V_\alpha\vDash \varphi(\gamma)\)
        lies below \(\eta\).\label{item:sc_correct}
        \item Every almost extendible cardinal is a limit of almost supercompact cardinals.\label{item:ex_to_limit_sc}
        \item If \(\eta\) is almost extendible,  \(\varphi\) is a formula, and \(\gamma < \eta\),
        then the least ordinal \(\alpha\) such that for all \(\beta \geq \alpha\), 
        \(V_\beta\vDash \varphi(\gamma,\alpha)\) lies below \(\eta\).
    \end{enumerate}
    \begin{proof}
        For \ref{item:ex_to_sc},
        fix \(\gamma < \eta < \nu\) and \(j : V_{\nu+1}\to V_{\nu'+1}\) such that \(j(\gamma) = \gamma\)
        and \(j(\eta) > \nu\).
        Then \(V_{\nu'+1}\) satisfies that for some \(\bar \nu < j(\eta)\),
        there is an elementary embedding
        \(\pi : V_{\bar \nu}\to V_{\nu'}\) such that \(\pi(\gamma) = \gamma\);
        namely, \(\pi = j\restriction V_{\bar \nu}\). Now by elementarity,
        for some \(\bar \nu < \eta\), there is an elementary embedding \(\pi : V_{\bar \nu}\to V_\nu\)
        such that \(\pi(\gamma) = \gamma\).

        \ref{item:sc_correct} is immediate: let \(\alpha'\) be an ordinal 
        such that \(V_{\alpha'}\vDash \varphi(\gamma)\). By almost supercompactness,
        there is some \(\alpha < \eta\) and an elementary \(\pi : V_\alpha\to V_\alpha'\)
        with \(\pi(\gamma) = \gamma\). In particular, by elementarity \(V_\alpha\vDash \varphi(\gamma)\).

        For \ref{item:ex_to_limit_sc}, it in fact suffices to assume that \(\eta\) is almost supercompact and 
        for all \(\gamma < \eta\),
        \(\eta\) is \((\gamma,\eta+1)\)-extendible.
        Fix \(\gamma < \eta\), and we will show there is an almost
        supercompact cardinal between \(\gamma\) and \(\eta\).
        Fix \(\pi : V_{\eta}\to V_{\eta'}\) such that \(\pi(\gamma) = \gamma\).
        Then in \(V_{\eta'}\), \(\eta\) is an almost supercompact cardinal between
        \(\gamma = \pi(\gamma)\) and \(\pi(\eta)\). So by elementarity,
        \(V_{\eta}\) satisfies that there is an almost supercompact
        cardinal \(\delta\) between \(\gamma\) and \(\eta\). Since \(\eta\) itself
        is almost supercompact, \ref{item:sc_correct} implies \(\delta\) really 
        is almost supercompact.
    \end{proof}
\end{lma}

\begin{lma}\label{lma:almost_super_crit}
    If \(\kappa\) is almost supercompact and \(\lambda\) is the least
    rank Berkeley cardinal, then for any ordinal \(\nu < \min(\kappa,\lambda)\),
    the almost supercompactness of \(\kappa\) is witnessed
    by elementary embeddings with critical point greater than \(\nu\).
    \begin{proof}
        Since \(\nu < \lambda\), there is an ordinal
        \(\alpha\) such that there is no elementary embedding
        \(j :V_\alpha\to V_\alpha\) with critical point less than \(\nu\).
        Since \(\kappa\) is almost supercompact, \cref{lma:extendible_to_sc} 
        implies that the least such ordinal \(\alpha\) is strictly less
        than \(\kappa\). Since the almost supercompactness of \(\kappa\) is witnessed
        by elementary embeddings fixing \(\alpha\),
        it is witnessed by elementary embeddings with critical point greater than \(\nu\).
    \end{proof}
\end{lma}

For one of our applications (the filter extension property, proved in \cref{section:filter_extension}),
we require a slight strengthening of these notions.
\begin{defn}\label{def:X_closed_rank_Berkeley} A cardinal \(\lambda\) is \textit{\(X\)-closed rank Berkeley}
if for all \(\gamma < \lambda < \alpha\), there is an elementary embedding
\(j :V_\alpha\to V_\alpha\) such that \(\gamma < \crit(j) < \lambda\)
and \(j(X) = j[X]\).\end{defn}
We note that if \(\lambda\) is \(X\)-closed rank Berkeley and
\(|Y|\leq^*|X|\), then \(\lambda\) is \(Y\)-closed rank Berkeley,
simply because if
\(f : X\to Y\) is a surjection, 
\(\alpha > \rank(X),\rank(Y)\) is an ordinal,
and \(j :V_\alpha\to V_\alpha\)
is an elementary embedding, then
\(j(Y) = j(f)[j(X)] = j(f)[j[X]] = j[Y]\).

If \(\kappa\) is a cardinal and \(X\in V_\kappa\),
then \(\kappa\) is \textit{\(X\)-closed almost extendible}
if for all \(\gamma < \kappa < \alpha\)
there is an elementary embedding \(j :V_\alpha\to V_{\alpha'}\)
such that \(j(\gamma) = \gamma\), \(j(X) = j[X]\), and \(j(\kappa) > \alpha\).

\begin{thm}\label{thm:closed_extendible}
    For any cardinal \(\lambda\), 
    for a closed unbounded class of \(\kappa\), \(\kappa\) is
    \(X\)-closed almost extendible for every \(X\) such that
    \(\lambda\) is \(X\)-closed rank Berkeley.
\end{thm}
The proofs above easily yield the following proposition:\footnote{
This is readily seen noting that given an embedding \(j : V_\alpha\to V_\alpha\),
the following are equivalent: (1) \(j(X) = j[X]\) and (2) letting \(\mathfrak{a} = |X|\) be the Scott cardinality
of \(X\),
\(j(\mathfrak a) = \mathfrak a\) and for all \(Y\in \mathfrak{a}\), \(j(Y) = j[Y]\).
In the proof of \cref{lma:lightface_to_extendible}, one will have to modify the definition of
\(\mathcal S\) by letting \(\eta_{\xi+1}\) be the least ordinal \(\nu > \eta_\xi\) such that
there is no \(\pi : V_{\nu}\to V_{\nu'}\) such that \(\pi(\eta_\xi) \geq \nu\),
\(\pi(\mathfrak a) = \mathfrak a\) and for all \(Y\in \mathfrak{a}\), \(\pi(Y) = \pi[Y]\).
It is then easy to see that if \(j : V_\alpha\to V_\alpha\) is elementary and \(j(\mathfrak{a}) = \mathfrak a\),
then \(j(\mathcal S) = \mathcal S\), which enables one to apply \cref{prp:bounded_vopenka}.
}
\begin{prp}\label{prp:X_closed}
    If there is an \(X\)-closed rank Berkeley cardinal,
    then there is a closed unbounded class of \(X\)-closed almost extendible cardinals.\qed
\end{prp}

\cref{thm:closed_extendible} is a trivial consequence of \cref{prp:X_closed} 
once one realizes that there is essentially just a set of \(X\) such that
\(\lambda\) is \(X\)-closed almost extendible.
We will use the following lemma, which will also be important later.
\begin{lma}\label{lma:j_closed}
    Suppose \(j : V_\alpha\to V_\alpha\) is an elementary embedding, 
    \(\kappa\) is almost supercompact, and \(\cf(\kappa) \geq \crit(j)\).
    Suppose \(A\in V_\alpha\)
    is a set such that \(j(A) = j[A]\).
    Then \(\scott(A) < \kappa\).
    \begin{proof}
        Let \(S\) be the set of Scott ranks of subsets of \(A\).
        Then \(j(S) = S\) since \(A\) and \(j(A)\) are in bijection.
        Moreover for all \(\nu\in S\), \(j(\nu) = \nu\)
        since for all \(B\subseteq A\), \(j(B) = j[B]\).
        Hence \(|S| < \crit(j)\).
        
        Let \(\xi = \sup(S\cap \kappa)\), and note that
        \(\xi < \kappa\) since \(\cf(\kappa) \geq \crit(j)\).

        Let \(\pi : V_{\bar \alpha}\to V_{\alpha}\)
        such that \(\xi < \bar \alpha < \kappa\),
        \(A\in \ran(\pi)\), and \(\pi(\xi) = \xi\).
        (See the comments following the definition of \((\gamma,\infty)\)-supercompactness.)
        Let \(\bar A = \pi^{-1}(A)\) 
        and let \(\nu\) be the Scott rank of \(\bar A\).
        Since \(\pi[\bar A]\subseteq A\),
        \(\nu \in S\). Therefore \(\nu \in S\cap \kappa\), which implies that \(\nu < \xi\),
        and hence \(\pi(\nu) < \pi(\xi) = \xi < \kappa\). 
        This completes the proof since \(\pi(\nu)\) is the Scott 
        rank of \(A\). 
    \end{proof}
\end{lma}
The assumption \(\cf(\kappa) \geq \crit(j)\) is not actually necessary, 
at least in the case of interest where \(\crit(j)\) is less than the least
rank Berkeley cardinal. The proof is slightly more complicated so we omit it,
but note that it allows us to remove this hypothesis for example from \cref{thm:antichains}.

\begin{proof}[Proof of \cref{thm:closed_extendible}]
    We may assume that \(\lambda\) is rank Berkeley,
    since otherwise the theorem is vacuous. 
    Applying \cref{cor:closed unbounded_extendible} and \cref{lma:extendible_to_sc}, 
    for each regular \(\gamma\), 
    let \(\rho_\gamma\) be the least almost supercompact cardinal of 
    cofinality \(\gamma\), and let \(\rho = \sup\{\rho_\gamma : \gamma\in \Reg\cap\lambda\}\).
    
    Let \(\Gamma\) denote the class of \(X\) such that \(\lambda\) is \(X\)-closed
    rank Berkeley.
    By \cref{lma:j_closed}, for each \(X\in \Gamma\), there is some \(Y\in V_\rho\)
    such that \(|X|\leq |Y|\). 
    If \(|X| \leq |Y|\) and \(\kappa\) is
    \(Y\)-closed almost extendible, then \(\kappa\) is \(X\)-closed almost extendible,
    so it suffices to show that there is a closed unbounded class of \(\kappa\) that
    is \(Y\)-closed almost extendible for all \(Y\in \Gamma\cap V_\rho\).
    This is an immediate consequence of \cref{prp:X_closed} and the closure of 
    closed unbounded classes under set-sized intersections.
\end{proof}
\subsection{Derived ultrafilters}
In this section, we show how to use the weak supercompactness
principles of the previous section
to derive the existence of certain kinds of ultrafilters.
It is assumed that the reader is familiar with the general
theory of normal fine ultrafilters
and supercompactness under the Axiom of Choice
(see for example \cite{Kanamori}); 
our exposition here is only intended to
highlight some minor modifications that are
needed to extend the theory to the choiceless setting.

Suppose \(\Lambda\) is a family of
sets and \(X = \bigcup_{\sigma\in \Lambda} \sigma\). 
A filter \(\mathcal F\) on \(\Lambda\)
is \textit{fine} if for all \(x\in \bigcup \Lambda\),
the set \(\{\sigma\in \Lambda : x\in \sigma\}\in \mathcal F\).
Thus an ultrafilter 
\(\mathcal U\) is fine if and only if \(j_\mathcal U[X]\subseteq \id_\mathcal U\)
where \(j_\mathcal U : V\to M_\mathcal U\) denotes the ultrapower of 
the universe by \(\mathcal U\).

The diagonal intersection of a family of sets \(\langle A_x : x\in X\rangle\) is 
the set \[\triangle_{x\in X}A_x = \left\{\sigma\subseteq P(X) : \textstyle\sigma\in \bigcap_{x\in \sigma} A_{x}\right\}\]
A filter \(\mathcal F\) on \(\Lambda\) is \textit{normal}
if it is closed under diagonal intersections: 
for any sequence \(\langle A_x : x\in X\rangle\subseteq \mathcal F\),
\(\triangle_{x\in X} A_x\in \mathcal F\). The Fodor
characterization of normality in terms of regressive functions,
familiar from set theory with choice, 
no longer suffices to characterize normality in ZF.
The usual argument instead shows the following:
\begin{lma}[Solovay]\label{lma:normal_regressive}
    Suppose \(\Lambda\) is a family of sets and
    \(X = \bigcup_{\sigma\in \Lambda} \sigma\). 
    If \(\mathcal F\) is a filter on \(\Lambda\),
    then \(\mathcal F\) is normal if and only if
    whenever \(T_\sigma\subseteq \sigma\)
    is nonempty for an \(\mathcal F\)-positive set of \(\sigma\in \Lambda\),
    there is some \(x\in X\) such that \(x\in T_\sigma\)
    for an \(\mathcal F\)-positive set of \(\sigma\in \Lambda\).
    \begin{proof}
        If \(\mathcal F\) is normal, then for each \(x\in X\),
        let \(A_x = \{\sigma \in \Lambda : x\notin T_\sigma\}\).
        Assume towards a contradiction that \(A_x\) is \(\mathcal F\)-null
        for all \(x\in X\). Then \(\triangle_{x\in X}A_x\in \mathcal F\) by normality,
        so fix some \(\sigma\in \triangle_{x\in X}A_x\). Then
        \(T_\sigma\subseteq \sigma\) is nonempty, but if \(y\in T_\sigma\),
        then \(y\in \sigma\), so \(\sigma\in A_y\) since \(\sigma\in \triangle_{x\in X}A_x\),
        but then \(y\notin T_\sigma\) by the definition of \(A_y\).

        Conversely, if \(\langle A_x : x\in X\rangle\subseteq \mathcal F\),
        let \(T_\sigma = \{x\in \sigma : \sigma\notin A_x\}\). If \(T_\sigma\neq\emptyset\)
        for an \(\mathcal F\)-positive set of \(\sigma\in \Lambda\), then there is some
        \(x\in X\) such that \(x\in T_\sigma\) for all \(\sigma\) in an \(\mathcal F\)-positive set \(S\subseteq \Lambda\).
        In other words, \(\sigma\notin A_x\) for all \(\sigma\in S\), or \(S\cap A_x = \emptyset\),
        which contradicts that \(A_x\in \mathcal F\) and \(S\) is \(\mathcal F\)-positive.
    \end{proof}
\end{lma}

A filter \(F\) is \textit{\(\gamma\)-descendingly closed}
if whenever \(\langle A_\alpha : \alpha < \gamma\rangle\)
is a decreasing sequence of sets in \(F\), \(\bigcap_{\alpha < \gamma} A_\alpha\in F\). 
A standard argument shows that if \(\eta\) is almost supercompact,
then for all \(\gamma < \eta\) and all sets \(Y\),
there is a \(\gamma\)-descendingly closed
normal fine ultrafilter on \(P(Y)\) 
that concentrates on the set of \(\sigma\subseteq Y\) such that
\(\theta(\sigma) < \eta\). 
\begin{defn}\label{def:derived}
    Suppose \(M\) and \(N\) are transitive sets, \(X\) is an element of \(M\),
    \(j : M\to N\) is an elementary embedding,
    and \(a\in j(X)\).
    The \textit{\(M\)-ultrafilter on \(X\) derived from \(j\) 
    using \(a\)} is \(\{A\in P(X)\cap M : a\in j(A)\}\).
\end{defn}
\begin{prp}\label{prp:normal_ultrafilters}
    Suppose \(\eta\) is \((\gamma,\infty)\)-supercompact
    and \(Y\) is a set. 
    \begin{enumerate}[(1)]
        \item \label{item:simple_uf}There is a \(\gamma\)-descendingly closed
    normal fine ultrafilter on \(P(Y)\) concentrating on the set
    of \(\sigma\subseteq Y\) such that \(\scott(\sigma) < \eta\).
        \item \label{item:ugly_uf} Set 
        \(\eta' = \eta\) if \(\eta\) is inaccessible and \(\eta' = \eta+1\) otherwise.
        Let \(\Lambda\) be the set of \(\sigma\subseteq Y\) such that 
        \(\dscott(\sigma) < \eta'\).
        Then there is a 
        \(\gamma\)-descendingly closed
        normal fine ultrafilter on \(P(\Lambda)\) concentrating on
        \(\sigma\subseteq P(\Lambda)\) such that \(\bigcup \sigma\in \Lambda\).
    \end{enumerate}
    \begin{proof}
    \ref{item:simple_uf}
    Fix a limit ordinal \(\nu\) such that
    \(Y\in V_\nu\). 
    By \((\gamma,\infty)\)-supercompactness, for some \(\bar \nu < \eta\),
    there is an elementary \(\pi : V_{\bar \nu}\to V_{\nu}\) 
    such that \(\pi(\gamma) = \gamma\)
    and \(\eta,Y\in \ran(\pi)\). (See \cref{foot:x_in_ran}.) 
    Let \(\bar Y = \pi^{-1}(Y)\) and \(\bar \eta = \pi^{-1}(\eta)\).
    Since \(|\pi[\bar Y]| = |\bar Y|\) and \(\bar Y\in V_\eta\), \(\scott(\pi[\bar Y]) < \eta\), 
    Let \(U\) be the \(\gamma\)-descendingly closed normal fine ultrafilter on 
    \(P(\bar Y)\) derived from \(\pi\) using \(\pi[\bar Y]\),
    so \(U\) concentrates on the set of \(\sigma\subseteq Y\) such that
    \(\scott(\sigma) < \bar \eta\). Then \(\pi(U)\)
    is a \(\gamma\)-descendingly closed
    normal fine ultrafilter concentrating on \(\sigma\subseteq Y\) such that \(\scott(\sigma) < \eta\).

    \ref{item:ugly_uf}
    Next, let \(\bar \Lambda = \pi^{-1}[\Lambda]\),
    and let \(W\) be the normal fine ultrafilter on 
    \(P(\bar \Lambda)\) derived from \(\pi\) using
    \(\pi[\bar \Lambda]\). We will show that \(\bigcup j[\bar \Lambda] \in \Lambda\),
    so that \(\pi(W)\) witnesses the final conclusion of the proposition.
    Define a partial surjection \(p : V_{\bar \nu}\times V_\eta\to \bigcup j[\bar \Lambda]\)
    by setting \(p(g,a) = \pi(g)(a)\) whenever \(g\) is a function with \(a\in \pi(\dom(g))\). 
    Then \(p\) witnesses that \(\dscott(\bigcup \sigma) \leq \eta\), 
    which implies \(\bigcup \sigma\in \Lambda\) unless \(\eta\) is inaccessible. 
    But if \(\eta\) is inaccessible, then
    \(\alpha = \sup_{\tau\in \bar \Lambda} \dscott(j(\tau))\) is less than \(\eta\), and 
    \(p[V_{\bar \nu} \times V_\alpha] = \bigcup j[\bar \Lambda]\). Therefore
    \(\dscott(\bigcup j[\bar \Lambda]) \leq \alpha < \eta\), as desired.
    \end{proof}
\end{prp}
\subsection{Regular cardinals}\label{section:regular}
If \((X,\preceq)\) is a wellfounded preorder, 
the \textit{rank of \(\preceq\)}, denoted by \(\rank({\preceq})\), 
is the least ordinal \(\alpha\)
admitting a function \(f : X\to \alpha\) 
such that \(x\preceq y\) implies \(x\leq y\).
For \(x\in X\), the \textit{rank of \(x\) in \(\preceq\)}, denoted by \(\rank_{{\preceq}}(x)\),
is the rank of the restriction of \(\preceq\) to \(\{y\in X : y\preceq x\text{ and }x\npreceq y\}\).

A prewellorder is a wellfounded preorder that is total; in other words,
it is a nonstrict wellorder without the antisymmetry condition. The prewellorders of a set
\(X\) are in one-to-one correspondence with the surjective functions from \(X\)
onto an ordinal.
If \(\Gamma\) is family of sets,
then \(\delta(\Gamma)\) denotes
the supremum of the ranks of all prewellorders in \(\Gamma\).
\begin{prp}\label{prp:delta_lemma}
    Suppose \(\Gamma\) is a set whose powerset carries 
    a \(\gamma\)-descendingly closed fine filter \(F\) that concentrates
    on the set of \(\sigma\subseteq\Gamma\) such that \(\delta(\sigma) < \delta(\Gamma)\).
    Then \(\cf(\delta(\Gamma)) \neq \gamma\).
    \begin{proof} 
        Let \(f : \gamma \to \delta(\Gamma)\) be an increasing function.
        We will show that \(f[\gamma]\) is bounded below \(\delta(\Gamma)\).
        For \(\alpha < \gamma\),
        let \(A_\alpha\) denote the set of all \(\sigma\in Y\)
        such that \(\rank(E) \geq f(\alpha)\) for some \(E\in \sigma\).
        Notice that \(A_\alpha \subseteq A_\beta\) for \(\beta\leq \alpha\)
        and by fineness, \(A_\alpha\in F\).
        Since \(F\) is \(\gamma\)-descendingly closed,
        \(A = \bigcap_{\alpha < \gamma} A_\alpha \in F\);
        note that \(A\) is the set of \(\sigma\subseteq \Gamma\) such that
        \(\delta(\sigma) \geq \sup f[\gamma]\).
        By our assumptions on \(F\), any \(F\)-large set contains some \(\sigma\) such that 
        \(\delta(\sigma) < \delta(\Gamma)\),
        so we may fix such a \(\sigma\) belongs to \(A\).
        Now \(\sup f[\gamma]\leq \delta(\sigma) < \delta(\Gamma)\),
        so \(f[\gamma]\) is bounded below \(\delta(\Gamma)\), as desired.
    \end{proof}
\end{prp}
Let \(\theta(X) = \delta(P(X\times X))\), the supremum of the
ranks of all prewellorders on \(X\).
Note that \(\theta(X)\) is the least ordinal 
that is not the range of a (partial) function on \(X\),
which is also known as the Lindenbaum number of \(X\).
\begin{cor}\label{cor:theta}
    Suppose that \(\gamma \leq \eta\) are cardinals, \(X\) is a set
    such that \(\eta \leq^* X\times X\leq^* X\). Assume
    there is a \(\gamma\)-descendingly closed fine filter on \(P(P(X))\)
    that concentrates on the set of \(\sigma\subseteq P(X)\) such that
    \(\theta(\sigma) < \eta\). Then \(\cf(\theta(X)) \neq \gamma\).
    \begin{proof}
        Since \(X\times X\leq^* X\), there is a \(\gamma\)-descendingly 
        closed fine filter \(\mathcal F\) on \(P(P(X\times X))\)
        that concentrates on the set of \(\sigma\subseteq P(X\times X)\) such that
        \(\theta(\sigma) < \eta\).
        Let \(\Gamma = P(X\times X)\).
        We will show that for all \(\sigma\subseteq \Gamma\), if
        \(\theta(\sigma) < \eta\), then
        \(\delta(\sigma) < \theta(X)\). Therefore the filter \(\mathcal F\)
        witnesses the hypothesis of \cref{prp:delta_lemma}, which implies the desired conclusion.

        Fix \(\sigma\subseteq \Gamma\) such that
        \(\theta(\sigma) < \eta\).
        For each \(x\in X\), let \(g_x : \sigma\to \delta(\sigma)\)
        be defined by \(g_x(E) = \rank_E(x)\). Let 
        \(A_x = g_x[\sigma]\) and let \(f_x : \alpha_x\to A_x\)
        be the increasing enumeration of \(A_x\). Since \(A_x\) is the surjective
        image of \(\sigma\),
        \(\alpha_x < \theta(\sigma) < \eta \leq \theta(X)\).
        Let \(g : X\to \theta(\sigma)\) be a surjection.
        Then define a partial surjection \(F : X\times X\to \delta(\sigma)\)
        by setting \(F(x,y) = f_x(g(y))\) whenever \(g(y) < \alpha_x\).
        It follows that \(\delta(\sigma) < \theta(X\times X) = \theta(X)\), as claimed.
    \end{proof}
\end{cor}

\begin{cor}
    If \(\eta\) is almost supercompact
    and \(X\) is a set such that \(\eta\leq^* X\times X\leq^* X\),
    then \(\cf(\theta(X)) \geq \eta\).\qed
\end{cor}
For any ordinal \(\gamma\), \(\gamma^+ = \theta(\gamma)\), so we conclude:
\begin{cor}\label{cor:successor}
    If \(\eta\) is almost supercompact,
    then every successor cardinal
    greater than or equal to \(\eta\)
    has cofinality at least \(\eta\).\qed
\end{cor}

\begin{cor}
    Assume there is an almost extendible cardinal.
    Then there is a proper class of regular cardinals.
    \begin{proof}
        Let \(\kappa\) be an almost extendible cardinal.
        By \cref{lma:extendible_to_sc},
        \(\kappa\) is a limit of almost supercompact cardinals,
        and by \cref{cor:successor}, if \(\eta\) is almost supercompact,
        the cofinality of \(\eta^+\) is either \(\eta\) or \(\eta^+\), which means
        that either \(\eta\) or \(\eta^+\) is a regular cardinal.
        Thus \(\kappa\) is a limit of regular cardinals. Now assume towards a contradiction
        that for some ordinal \(\alpha\), there are no regular cardinals greater than \(\alpha\).
        By \cref{lma:extendible_to_sc}, the least such \(\alpha\) is less than \(\kappa\),
        but this contradicts that \(\kappa\) is a limit of regular cardinals.
    \end{proof}
\end{cor}
Combining this with \cref{cor:closed unbounded_extendible}, we obtain the answer to Asper\'o's question:
\begin{thm}\label{thm:regulars}
    Assume there is a rank Berkeley cardinal. Then there is a proper class of regular cardinals.\qed
\end{thm}
\subsection{The wellordered collection lemma}
In this section, we show that large cardinal axioms imply weak versions of the Axiom of Choice for wellordered
families. The main combinatorial result is the following. The proof is reminiscent of
Woodin's Coding Lemma for \(L(V_{\lambda+1})\) \cite{SEM2}.
\begin{thm}\label{thm:wo_collection_filter}
    Suppose \(\Lambda\) is a family of sets,
    \(X = \bigcup \Lambda\), and \(\beta\) is an ordinal
    such that for all \(\alpha\leq \beta\),
    there is an \(\alpha\)-descendingly closed fine filter \(F\) on \(P(\Lambda)\)
    concentrating on the set of \(\sigma\subseteq \Lambda\) such that \(\bigcup \sigma\in \Lambda\).
    Then for any sequence \(\langle S_\xi : \xi < \beta\rangle\)
    of nonempty subsets of \(X\), there is a set \(A\in \Lambda\) such that
    \(S_\xi\cap A\neq \emptyset\)
    for all \(\xi < \beta\).
    \begin{proof}
        By induction, we will show that for any ordinal \(\alpha \leq \beta\),
        there is a set \(A\in \Lambda\) such that
        \(S_\xi\cap A\neq \emptyset\)
        for all \(\xi < \alpha\).

        Assume first that the induction hypothesis holds for some \(\alpha < \beta\),
        and we will prove it for \(\alpha+1\).
        Fix a set \(A\in \Lambda\) such that \(S_\xi\cap A\neq\emptyset\) for all
        \(\xi < \alpha\).
        Since \(X = \bigcup \Lambda\),
        there is a set \(B\in \Lambda\) such that
        \(x\in B\). Since there is a fine filter on
        \(\Lambda\) concentrating on the set of \(\sigma\subseteq \Lambda\) such that \(\bigcup \sigma\in \Lambda\),
        there is some \(C\in \Lambda\) such that \(A\cup B\subseteq C\) (since almost all \(\sigma\subseteq \Lambda\)
        have \(A,B\in \sigma\)).
        It follows that
        for \(F\)-almost all \(\sigma\subseteq \Lambda\), 
        \(S_\alpha\cap C\neq \emptyset\). Therefore
        \(F\) witnesses the induction hypothesis for \(\alpha+1\).
        
        Assume instead that \(\gamma \leq \beta\) is a limit ordinal, 
        and the induction hypothesis holds for all \(\alpha < \gamma\),
        and we will prove the induction hypothesis for \(\gamma\). 
        Let \(F\) be a \(\gamma\)-descendingly closed fine filter on \(P(\Lambda)\)
        concentrating on the set of \(\sigma\subseteq \Lambda\) such that \(\bigcup \sigma\in \Lambda\).
        For \(\alpha < \gamma\), let \(\mathcal A_\alpha\subseteq P(\Lambda)\) denote the set of 
        \(\sigma\subseteq \Lambda\) such that for all \(\xi < \alpha\),
        \(S_\xi\cap \bigcup \sigma\neq \emptyset\). Clearly the sequence 
        \(\langle \mathcal A_\alpha : \alpha < \gamma\rangle\) is descending.
        For each \(\alpha < \gamma\), our induction hypothesis implies that
        there is a set \(A\in \Lambda\) such that \(S_\xi \cap A\neq \emptyset\) for all
        \(\xi < \alpha\). By fineness, \(F\) concentrates on the set of \(\sigma\subseteq P(\Lambda)\)
        such that \(A\in \sigma\), and so \(\mathcal A_\alpha\in F\). 
        Since \(F\) is \(\gamma\)-descendingly closed, 
        \(\mathcal A = \bigcap_{\alpha < \gamma} \mathcal A_\alpha\) belongs to \(F\),
        or in other words,
        \(F\) concentrates on the set of \(\sigma\subseteq \Lambda\) such that 
        \(S_\xi\cap (\bigcup \sigma)\neq \emptyset\) for all \(\xi < \gamma\).
        Since \(F\) is a (proper) filter, there is some such \(\sigma\) satisfying 
        \(\bigcup \sigma\in \Lambda\),
        and so \(A = \bigcup \sigma\) witnesses the induction hypothesis for \(\gamma\). 
    \end{proof}
\end{thm}

\begin{cor}[Wellordered collection lemma]\label{thm:wo_collection}
    Suppose \(\eta\) is 
    a cardinal that is 
    \((\gamma,\infty)\)-supercompact for all regular \(\gamma \leq \beta\).
    Then for any sequence \(\langle S_\xi : \xi < \beta\rangle\) of nonempty sets,
    there is a set \(\tau\) such that \(\dscott(\tau) \leq \eta\)
    and \(S_\xi\cap \tau\neq \emptyset\) for all \(\xi < \beta\).
    If \(\eta\) is inaccessible, one can find such a \(\tau\) with \(\scott(\tau) < \eta\).
    \begin{proof}
        Let \(X = \bigcup_{\xi < \beta} S_\xi\).
        Let \(\Lambda\) be the set of \(\tau\subseteq X\) such that \(\dscott(\tau) < \eta'\).
        For any regular cardinal \(\gamma \leq \beta\), 
        \cref{prp:normal_ultrafilters} \ref{item:ugly_uf} 
        yields a \(\gamma\)-descendingly closed
        normal fine ultrafilter on \(P(\Lambda)\) concentrating on the set of
        \(\sigma\subseteq \Lambda\) such that \(\bigcup \sigma\in \Lambda\). 
        Applying \cref{thm:wo_collection_filter} then yields the theorem.
        (Note here that \(\dscott(\tau) \leq\scott(\tau)+1\), so if \(\eta\)
        is a limit ordinal, \(\dscott(\tau) < \eta\) implies \(\scott(\tau) <\eta\).)
        \end{proof}
\end{cor}

\begin{prp}\label{thm:wo_collection_extendible}
    Suppose \(\kappa\) is almost extendible. 
    Then for any \(\beta < \kappa\), 
    for any sequence \(\langle S_\xi : \xi < \beta\rangle\),
    there is a set \(\sigma\) with \(\scott(\sigma) < \kappa\)
    such that \(S_\xi\cap \sigma\neq\emptyset\) for all \(\xi < \beta\).
    \begin{proof}
        This follows from the fact that almost extendible cardinals
        are limits of almost supercompact cardinals (\cref{lma:extendible_to_sc}).
    \end{proof}
\end{prp}

\begin{thm}\label{thm:wellordered_collection}
    Suppose there is a rank Berkeley cardinal.
    For every ordinal \(\kappa\), there is a Scott
    cardinal \(\mathfrak{a}\) such that
    for any sequence \(\langle A_\alpha : \alpha < \kappa\rangle\) of nonempty
    sets, there is a set \(\sigma\) of cardinality at most \(\mathfrak{a}\) and 
    \(A_\alpha\cap \sigma\neq \emptyset\) for all \(\alpha < \kappa\).\qed
\end{thm}
\section{Measurable cardinals}
\subsection{Filter bases}
A \textit{filter base} is a family of nonempty sets \(\mathcal B\)
that is downwards directed under inclusion.
If \(X\) is a set,
a filter base \(\mathcal B\) is \textit{\(X\)-closed}
if for any \(\langle A_x : x \in X\rangle\subseteq \mathcal B\),
there is some \(A\in \mathcal B\) such that \(A\subseteq \bigcap_{x \in X} A_x\).
If \(Y\) is a family of sets, \(\mathcal B\) is \textit{\(Y\)-complete}
if it is \(X\)-closed for all \(X\in Y\).
A filter base \(F\) 
is a \textit{filter} if it is closed upwards under inclusion.
The \textit{filter generated by a filter base \(\mathcal B\)}
is the family of sets that contain some element of \(\mathcal B\). 

If \(F\) is a filter on \(X\), a set \(A\) is
\textit{\(F\)-null} if \(X\setminus A\) belongs to \(F\). 
The dual ideal of \(F\), denoted by \(F^*\),
is the set of \(A\subseteq X\) such that \(X\setminus A\in F\).
A set \(S\) is \textit{\(F\)-positive} if \(S\cap A\neq \emptyset\)
for all \(A\in F\), or equivalently, if \(S\) is not \(F\)-null.
The set of \(F\)-positive subsets of \(X\) is denoted by \(F^+\).
If \(F\) is a filter and \(S\in F^+\), \(F\restriction S\)
denotes the filter generated by \(F\cup \{S\}\). 

The following theorem on the completeness
of the filter generated by a filter base is almost
a restatement of \cref{thm:wo_collection_filter}.
\begin{thm}\label{thm:filter_base}
    Suppose \(\mathcal B\) is a filter base on a set \(X\) and \(\beta\) is an ordinal
    such that for all regular \(\gamma \leq \beta\),
    there is a \(\gamma\)-descendingly closed fine filter on \(P(\mathcal B)\)
    concentrating on the set of \(\sigma\subseteq \mathcal B\) such that \(\bigcap \sigma\in \mathcal B\).
    Then the filter generated by \(\mathcal B\) is \(\beta\)-closed.
    \begin{proof}
        Let \(F\) be the filter generated by \(\mathcal B\), and suppose \(\langle A_\xi : \xi < \beta\rangle\)
        is a set of elements of \(F\). We will find a set \(A\in \mathcal B\)
        such that \(A\subseteq \bigcap_{\xi < \beta} A_\xi\). 
        Let \(S_\xi = \{A\in \mathcal B : A\subseteq A_\xi\}\)
        so that \(\langle S_\xi : \xi <\beta\rangle\) is a sequence of nonempty subsets of \(\mathcal B\).
        Let \(\Lambda\) be the set of \(\sigma\subseteq \mathcal B\) such that \(\bigcap \sigma\in \mathcal B\).
        Then by \cref{thm:wo_collection_filter}, there is some \(\sigma\in \Lambda\)
        such that \(S_\xi\cap \sigma\neq \emptyset\) for all \(\xi < \beta\). Letting
        \(A = \bigcap \sigma\), since \(\sigma\in \Lambda\), \(A\in \mathcal B\), and
        since \(\sigma\cap S_\xi\neq \emptyset\) for all \(\xi < \beta\), \(A\subseteq \bigcap_{\xi < \beta}A_\xi\).
    \end{proof}
\end{thm}
\begin{cor}\label{thm:filter_base_sc}
    Suppose \(\beta < \eta\) are ordinals
    such that \(\eta\) is \((\gamma,\infty)\)-supercompact for all regular \(\gamma \leq \beta\).
    Then any \(V_\eta\)-complete filter base generates a \(\beta\)-closed filter.\qed
\end{cor}
\begin{thm}\label{thm:completeness}
    Suppose \(\eta\) is almost supercompact.
    Then for all ordinals \(\delta\) of cofinality at least \(\eta\),
    the closed unbounded filter on \(\delta\) is \(\eta\)-complete.
    \begin{proof}
        By \cref{thm:filter_base_sc}, 
        it suffices to show that the set of closed unbounded subsets of \(\delta\) is 
        a \(V_\eta\)-complete filter base. The proof of this fact is
        familiar from the standard theory of closed unbounded sets.
        Suppose \(X\in V_\eta\) and \(\langle C_x : x\in X\rangle\) 
        is a sequence of closed unbounded subsets of \(\delta\).
        The set \(\bigcap_{x\in X} C_x\) is clearly closed,
        so it suffices to show it is unbounded.

        Fix \(\alpha_0 < \delta\), and we will exhibit an ordinal \(\alpha_\omega > \alpha_0\)
        that belongs to \(\bigcap_{x\in X} C_x\). For \(n < \omega\) and \(x\in X\),
        let \(\alpha_{n+1}(x)\) be the least element of \(C_x\) above \(\alpha_n\),
        and let \(\alpha_{n+1} = \sup_{x\in X}\alpha_{n+1}(x)\).
        Since \(\cf(\delta) \geq \eta > \theta(X)\), 
        \(\alpha_{n}\) is defined for all \(n < \omega\).
        Let \(\alpha_\omega = \sup_{n < \omega}\alpha_n\).
        Then \(\alpha_\omega\) is a limit point of \(C_x\) for all \(x\in X\),
        and therefore \(\alpha_\omega\in \bigcap_{x\in X} C_x\).
    \end{proof}
\end{thm}

\subsection{Wellfounded filters}
A filter \(F\) on a set \(X\) is \textit{\(\gamma\)-wellfounded} if the 
reduced product \(\gamma^X/F\) is wellfounded.
We say that \(F\) is \textit{wellfounded} if it is \(\gamma\)-wellfounded
for all ordinals \(\gamma\).
If \(F\) is a \(\gamma\)-wellfounded filter on \(X\)
and \(\langle \mathcal M_x : x\in X\rangle\) is a sequence of wellfounded 
structures of rank at most \(\gamma\),
then the reduced product \(\mathcal M = \prod_{x\in X}\mathcal M_x / F\) is again wellfounded: 
define \(o : \prod_{x\in X}\mathcal M_x\to \gamma^X\)
by setting \(o(f)(x) = \rank_{\mathcal M_x}(f(x))\), and note that \(o\) induces
a rank function \(o: \mathcal M\to \gamma^X/F\). 

A set \(A\subseteq P(X)\) lies below a set \(B\subseteq P(Y)\)
in the \textit{Kat\v{e}tov order}, denoted \(A\kat B\), if 
there is a function \(f : Y\to X\) such that for all \(S\in A\), \(f^{-1}[S]\in B\).
If \(F\kat G\) where \(G\) is a \(\gamma\)-wellfounded filter, then \(F\) is
a \(\gamma\)-wellfounded filter as well.

The following universality fact for fine filters, due to Kunen,
allows us to conclude that under large cardinal hypotheses
all sufficiently complete filters are wellfounded.
If \(\mathcal B\) is a filter base on \(X\) and \(\sigma\) is a subset of \(P(X)\),
let \(A_{\mathcal B}(\sigma) = \bigcap_{A\in \mathcal B\cap \sigma} A\).
If \(X\) is wellordered, define a partial function \(\chi_{\mathcal B} : P(P(X)) \to X\)
by \(\chi_{\mathcal B}(\sigma) = \min(A_{\mathcal B}(\sigma)))\).
\begin{thm}[Kunen]\label{thm:katetov}
    Assume \(\delta\) is an ordinal,
    \(\mathcal B\) is a filter base on \(\delta\), 
    and \(\mathcal W\) is a fine filter on \(P(P(\delta))\)
    concentrating on the set \(\Gamma\) of all \(\sigma\subseteq P(\delta)\) such that
    \(\bigcap (\mathcal B\cap \sigma)\neq \emptyset\). 
    Then \(\mathcal B\) lies below \(\mathcal W\)
    in the Kat\v{e}tov order.
    \begin{proof}
        Note that \(\chi_\mathcal B\) is defined on the \(\mathcal W\)-large set \(\Gamma\), and
        for all \(A\in \mathcal B\), 
        \(\chi_{\mathcal B}^{-1}[A] \in \mathcal W\)
        since \(\mathcal W\) is fine and
        \(\{\sigma\in P(P(\delta)) : A\in \sigma\}\subseteq \chi_{\mathcal B}^{-1}[A]\).
    \end{proof}
\end{thm}

The following lemma on the completeness of filters on ordinals
is useful to keep in mind.
\begin{lma}\label{lma:ordinal_complete}
    Suppose \(\kappa\) is an ordinal, \(X\) is a set,
    and there is no \(\kappa\)-sequence of distinct
    subsets of \(X\).
    If \(F\) is a \(\kappa\)-complete
    filter on an ordinal \(\delta\), then 
    \(F\) is \(X\)-closed.
    \begin{proof}
        Suppose \(\langle S_x:x\in X\rangle\) is a sequence
        of sets such
        that \(\bigcup_{x\in X} S_x\) is \(F\)-positive,
        and we will show that \(S_x\) is \(F\)-positive for some \(x\).

        Let \(S = \bigcup_{x\in X} S_x\).
        For each \(\alpha \in S\),
        \(D_\alpha = \{x\in X : \alpha\in S_x\}\).
        Then \(\{D_\alpha : \alpha < \delta\}\) is a wellorderable family of
        subsets of \(X\), and hence it has cardinality less than \(\kappa\).

        Let \(A_\alpha = \{\xi < \delta : D_\xi = D_\alpha\}\). 
        Then \(|\{A_\alpha : \alpha \in S\}| < \kappa\) and
        \(\bigcup_{\alpha\in S} A_\alpha = S\) since \(\alpha\in A_\alpha\).
        Thus \(\{A_\alpha : \alpha\in S\}\) is a
        partition of \(S\) into fewer than \(\kappa\)-many sets, 
        and so since \(F\) is \(\kappa\)-complete
        there is some \(\alpha\in S\) such that \(A_\alpha\) is \(F\)-positive.
        Fix \(x\in X\) such that \(\alpha\in S_x\), and note that
        \(A_\alpha\subseteq S_x\), and hence \(S_x\) is \(F\)-positive, as desired.
    \end{proof}
\end{lma}
If \(X\) is a set, the \textit{Hartogs number of \(X\),} denoted \(\aleph(X)\), is the 
least ordinal \(\alpha\) such that there is no injection from \(\alpha\) to \(X\).
\begin{lma}\label{lma:complete_wf}
    Suppose there is a wellfounded \(\nu\)-complete
    fine ultrafilter \(\mathcal W\) on \(P(P(\delta))\)
    that concentrates on the set of \(\sigma\subseteq P(\delta)\) such that
    \(\aleph(P(\sigma)) < \epsilon\).
    Then every \(\epsilon\)-complete filter \(F\) on \(\delta\)
    extends to a wellfounded \(\nu\)-complete ultrafilter.
    Moreover the set of \(\epsilon\)-complete ultrafilters on \(\delta\)
    can be wellordered.
    \begin{proof}
        \cref{lma:ordinal_complete} puts us in a position to apply
        \cref{thm:katetov}. The point is that if \(F\) is an \(\epsilon\)-complete filter on \(\delta\)
        and \(\sigma\subseteq P(\delta)\) has the property that \(\aleph(P(\sigma)) < \epsilon\),
        then there is no \(\epsilon\)-sequence of distinct subsets of \(\sigma\). Applying
        \cref{lma:ordinal_complete} with \(\kappa = \epsilon\), we have that \(F\) is 
        \(\sigma\)-closed. Therefore since \(\mathcal W\) concentrates on the set of 
        \(\sigma\subseteq P(\delta)\) such that
        \(\aleph(P(\sigma)) < \epsilon\), in fact \(\mathcal W\) concentrates on 
        the set of \(\sigma\subseteq P(\delta)\) such that \(\bigcap (\sigma\cap F)\in F\),
        and so by \cref{thm:katetov},
        \(F\) lies below \(\mathcal W\) in the Kat\v{e}tov order.
        Let \(f : P(P(\delta))\to \delta\) be such that
        for all \(A\in F\), \(f^{-1}[A]\in \mathcal W\).
        Then the ultrafilter \(\{A\subseteq \delta : f^{-1}[A]\in \mathcal W\}\)
        extends \(F\).

        Define a binary relation \(\leq\) on the set of \(\epsilon\)-complete ultrafilters on \(\delta\)
        by setting \(U_0 \leq U_1\) if \(\chi_{U_0}(\sigma) \leq \chi_{U_1}(\sigma)\) for
        \(\mathcal W\)-almost all \(\sigma\). (See the paragraph before \cref{thm:katetov} for the
        notation \(\chi_U\).) 
        This relation is transitive because \(\mathcal W\) is a filter,
        antisymmetric because \(\mathcal W\) is fine, 
        linear because \(\mathcal W\) is maximal, 
        and wellfounded because \(\mathcal W\) is wellfounded.
    \end{proof}
\end{lma}
The hypothesis of \cref{lma:complete_wf} is a combinatorial consequence of
the existence of a \((0,\infty)\)-supercompact cardinal
whose associated elementary embeddings
have critical point greater than \(\nu\). For the most part,
this follows from \cref{prp:normal_ultrafilters}. The only
new element is showing that one can obtain a \textit{wellfounded} ultrafilter,
which follows from the proof of \cref{prp:normal_ultrafilters} combined
with the following lemma:
\begin{lma}
    Suppose \(\nu\) is an ordinal such that \(V_\nu\preceq_{\Sigma_1} V\)
    and \(\pi : V_{\bar \nu}\to V_\nu\) is an elementary embedding.
    If \(D\) is an ultrafilter derived from \(\pi\),
    then \(\pi(D)\) is wellfounded.
    \begin{proof}
        Suppose \(X\in V_{\bar \nu}\), \(a\in \pi(X)\), and \(D = \{A\subseteq X : a\in \pi(A)\}\)
        is the ultrafilter on \(X\) derived from \(\pi\) using \(A\).
        Then \(D\) is \(\alpha\)-wellfounded for all \(\alpha < \bar \nu\)
        because the ultrapower of \(\alpha\) by \(D\) can be embedded into \(\nu\),
        by setting \([f]_D = \pi(f)(a)\) for any \(f : X\to \alpha\). 
        By elementarity, \(\pi(D)\) is \(\alpha\)-wellfounded for all \(\alpha < \nu\).
        It follows that \(\pi(D)\) is wellfounded in \(V_\nu\), and so since \(V_\nu\preceq_{\Sigma_1} V\),
        \(\pi(D)\) is wellfounded.
    \end{proof}
\end{lma}

The proof of \cref{lma:complete_wf} derives from Kunen's theorem on ultrafilters under AD:
\begin{thm*}[Kunen]
    Assume \(\AD + \textnormal{DC}_\mathbb R\).
    Then the set of all ultrafilters on 
    ordinals less than \(\theta(\mathbb R)\) is wellorderable.
    In fact, every ultrafilter on an ordinal less than \(\theta(\mathbb R)\) is
    ordinal definable.\qed
\end{thm*}
In fact, we will soon show that rank Berkeley cardinals
almost imply the ordinal definability of sufficiently complete ultrafilters (\cref{KOSection}).
\subsection{Saturated filters in normal measures}
A theorem of Woodin states that if \(j : V\to V\) is an elementary embedding
and \(\lambda\)-DC holds where \(\lambda = \kappa_\omega(j)\), 
then \(\lambda^+\) is measurable.
One can show without any choice assumptions that the existence of \(j:V\to V\)
implies that there can be no partition of
\(S^{\lambda^+}_\omega = \{\alpha < \lambda^+ : \cf(\alpha) = \omega\}\) into \(\lambda\)-many stationary sets.
Assuming \(\lambda\)-DC, one can then apply Ulam's splitting theorem: suppose \(2^{<\lambda} = \lambda\) and
\(F\) is a
\(\lambda^+\)-complete filter on \(X\) such that there is no
partition of \(X\) into \(\lambda\)-many disjoint \(F\)-positive sets; then there is a
partition \(\langle A_\xi : \xi < \nu\rangle\) of \(X\) into \(F\)-positive
sets such that for all \(\xi < \nu\),
\(F\restriction A_\xi\) is an ultrafilter. 
In the case that and \(F\) is the \(\omega\)-closed unbounded filter on \(\lambda^+\), it follows 
that there is a stationary set \(A\subseteq S^{\lambda^+}_\omega\)
such that \(F\restriction A\) is a \(\lambda^+\)-complete ultrafilter on \(\lambda^+\).
(The \(\lambda^+\)-completeness of the closed unbounded filter follows
from \(\lambda\)-DC, or even just \(\lambda\)-AC.)

Woodin's non-splitting lemma seems very specific to the closed unbounded filter, 
but here we will show that it actually applies much more generally (\cref{cor:weak_saturation}).
Ulam's argument seems to use its dependent choice hypothesis too heavily to be of any use
in the context of bare ZF,
but we will show that it \textit{can} be carried out
using almost supercompactness instead of DC. 
Then in \cref{section:closed unboundeds}, we will combine these results to prove the existence of 
a proper class of measurable
cardinals assuming a rank Berkeley cardinal.

If \(\mathcal E\) is a family of sets, a set \(S\) is
\textit{\(\mathcal E\)-positive} if \(S\cap A\neq \emptyset\) for all \(A\in \mathcal E\). 
\begin{lma}\label{lma:saturation}
    Suppose \(\mathcal U\) is a normal fine ultrafilter
    on \(P(X)\) and \(\mathcal E\in \mathcal U\).
    Suppose there is an \(I\)-indexed family
    of disjoint \(\mathcal E\)-positive sets.
    Then every function \(f : P(X)\to I\) is constant
    on a set in \(\mathcal U\).
    \begin{proof}
        Let \(f : P(X)\to I\) be a function.
        Suppose \(\langle S_i : i \in I\rangle\) is a family
        of \(\mathcal E\)-positive sets. 
        For \(\mathcal U\)-almost all \(\sigma\),
        \(\sigma\cap S_{f(\sigma)} \neq \emptyset\),
        and so applying \cref{lma:normal_regressive}, 
        since \(\mathcal U\) is normal,
        there is some \(x\in X\) such that
        for \(\mathcal U\)-almost all \(\sigma\),
        \(x\in S_{f(\sigma)}\). Since the sets 
        \(\langle S_i : i \in I\rangle\) are disjoint,
        there is a unique \(i\in I\) 
        such that \(x\in S_i\).
        Thus for \(\mathcal U\)-almost all \(\sigma\),
        \(f(\sigma) = i\). 
    \end{proof}
\end{lma}
For any ultrafilter \(\mathcal U\), let \(\kappa_\mathcal U\) denote the
largest cardinal \(\kappa\) such that \(\mathcal U\) is \(\kappa\)-complete.
If \(\mathcal U\) is a normal fine ultrafilter on \(P(\delta)\), either \(\kappa_\mathcal U\)
is defined or \(\mathcal U\) is principal.

A family of sets \(\mathcal E\) is \textit{\(\gamma\)-weakly saturated}
if there is no collection of disjoint \(\mathcal E\)-positive sets of cardinality \(\gamma\).
In our main applications, we will have \(\mathcal E\subseteq P(\delta)\) for some ordinal \(\delta\),
so that \(\mathcal E\) is wellorderable.
\begin{cor}\label{cor:weak_saturation}
    Suppose \(\delta\) is an ordinal, \(\mathcal U\) is a nonprincipal
    normal fine ultrafilter
    on \(P(\delta)\), and \(\mathcal E\in \mathcal U\). 
    Then \(\mathcal E\) is \(\kappa_\mathcal U\)-weakly saturated.
    \begin{proof}
        If \(P\) is a family of disjoint \(\mathcal E\)-positive sets,
        then \(P\) is wellorderable since there is an injection
        \(f:P\to\delta\) defined by \(f(S) = \min(S\cap \delta)\).
        By \cref{lma:saturation}, no family of disjoint \(\mathcal E\)-positive
        sets has cardinality \(\kappa_\mathcal U\), and therefore
        every family of disjoint \(\mathcal E\)-positive sets
        has cardinality less than \(\kappa_\mathcal U\). 
    \end{proof}
\end{cor}
A filter \(F\) is \textit{\((\kappa,\epsilon)\)-indecomposable}
if every family \(P\) of disjoint sets with \(|P| < \epsilon\) has a subfamily
\(Q\) such that \(|Q| < \kappa\) and \(\bigcup (P\setminus Q)\) is \(F\)-null.
Note that if \(F\) is \((\kappa,\epsilon)\)-indecomposable, so is every extension of \(F\).
\begin{lma}
    Suppose \(F\) is a filter that is \(\epsilon\)-complete and \(\kappa\)-weakly saturated.
    Then \(F\) is \((\kappa,\epsilon)\)-indecomposable.
    \begin{proof}
        Suppose \(P\) is a family of disjoint sets with \(|P| < \epsilon\).
        Let \(Q\subseteq P\) be the collection of \(F\)-positive sets in \(P\).
        Since \(F\) is \(\kappa\)-weakly saturated, \(|Q| < \kappa\).
        Since \(F\) is \(\epsilon\)-complete, \(\bigcup (P\setminus Q)\) is \(F\)-null,
        being the union of fewer than \(\epsilon\)-many 
        \(F\)-null sets.
    \end{proof}
\end{lma}

\begin{cor}\label{cor:extension_complete}
    Suppose \(\delta\) is an ordinal and \(\mathcal U\) is a nonprincipal
    normal fine ultrafilter
    on \(P(\delta)\). Suppose \(F\) is an \(\epsilon\)-complete
    filter in \(\mathcal U\). Then every \(\kappa_\mathcal U\)-complete
    filter extending \(F\) is \(\epsilon\)-complete.
    \begin{proof}
        Suppose \(G\) is a \(\kappa_\mathcal U\)-complete filter extending \(F\).
        Since \(F\in\mathcal U\),
        \(F\) is \(\kappa_\mathcal U\)-weakly saturated,
        and therefore \(F\) is \((\kappa_\mathcal U,\epsilon)\)-indecomposable,
        and hence so is \(G\). Since \(G\) is \((\kappa_\mathcal U,\epsilon)\)-indecomposable
        and \(\kappa_\mathcal U\)-complete, \(G\) is \(\epsilon\)-complete.
    \end{proof}
\end{cor}
If \(F\) is a filter, an \(F\)-positive set \(A\) is an \textit{atom of \(F\)}
if \(F\restriction A\) is an ultrafilter. In other words,
\(A\) cannot be partitioned into distinct \(F\)-positive sets.

\begin{thm}\label{thm:saturation}
    Suppose \(\delta\) and \(\epsilon\) are ordinals, \(\mathcal U\) is a nonprincipal
    normal fine ultrafilter on \(P(\delta)\), and \(F\in \mathcal U\) is an
    \(\epsilon\)-complete filter.
    Suppose there is a wellfounded \(\kappa_\mathcal U^+\)-complete
    fine ultrafilter \(\mathcal W\) on \(P(P(\delta))\)
    that concentrates on the set of \(\sigma\subseteq P(\delta)\) such that
    \(\aleph(P(\sigma)) < \epsilon\).
    Then \(\delta\) can be partitioned 
    into fewer than \(\kappa_\mathcal U\)-many atoms of \(F\).
    \begin{proof}
        Let \(\mathcal S\) be the set of \(\epsilon\)-complete
        ultrafilters extending \(F\).
        We start by showing 
        that for any \(\mathcal T\subseteq
        \mathcal S\) such that \(|\mathcal T| \leq \kappa_\mathcal U\), 
        there is a family
        \(\mathcal P\) of disjoint subsets of \(\delta\) such that each
        \(|\mathcal P\cap U| = 1\) for all \(U\in \mathcal T\).

        Since \(\mathcal W\) is fine, for any \(U_0,U_1\in \mathcal S\), 
        for \(\mathcal W\)-almost all \(\sigma\),
        there is some \(A\in \sigma\) such that \(A\in U_0\setminus U_1\).
        Since \(\mathcal W\) is \(\kappa_\mathcal U^+\)-complete and
        \(|\mathcal T| \leq \kappa_\mathcal U\), the quantifiers can be exchanged:
        that is, for \(\mathcal W\)-almost all \(\sigma\),
        for all \(U_0,U_1\in \mathcal T\), 
        there is an \(A\in \sigma\) such that \(A\in U_0\setminus U_1\).
        Fix such a set \(\sigma\) such that \(\aleph(P(\sigma)) < \epsilon\).
        Then \(\mathcal P = \{A_U(\sigma) : U\in \mathcal T\}\) is as desired.

        By \cref{lma:complete_wf}, \(\mathcal S\) can be wellordered.
        It follows that \(|\mathcal S| < \kappa_\mathcal U\).
        Assume not. Then there is some \(\mathcal T\subseteq\mathcal S\)
        such that \(|\mathcal T| = \kappa_\mathcal U\). Fix
        a family \(\mathcal P\) of disjoint subsets of \(\delta\) such that 
        \(|\mathcal P\cap U| = 1\) for all \(U\in \mathcal T\).
        Then \(\mathcal P\) is a family of \(\kappa_\mathcal U\)-many disjoint 
        \(F\)-positive sets, the existence of which contradicts \cref{lma:saturation}.

        Thus \(|\mathcal S| < \kappa_\mathcal U\), and so there is a family
        \(\mathcal P\) of disjoint subsets of \(\delta\) such that
        \(|\mathcal P\cap U| = 1\) for all \(U\in \mathcal S\). 
        We claim that \(\mathcal P\) partitions an \(F\)-large set into atoms of \(F\),
        which easily implies the conclusion of the lemma.

        We first show that \(\bigcup_{A\in \mathcal P} A\in F\). 
        Suppose not.
        Let \(S = \delta\setminus \bigcup_{A\in \mathcal P} A\).
        Then \(F\restriction S\) is an \(\epsilon\)-complete filter and
        so as a consequence of \cref{lma:complete_wf} 
        extends to a \(\kappa_\mathcal U^+\)-complete ultrafilter \(U\).
        By \cref{cor:extension_complete}, \(U\) is
        \(\epsilon\)-complete. Therefore for some
        \(A\in \mathcal P\), \(A\in U\). Since \(U\) is a proper filter, 
        \(S\cap A\neq \emptyset\), which contradicts that 
        \(S\subseteq \delta\setminus A\).

        By a similar argument, we show that each \(A\in \mathcal P\) is an \(F\)-atom.
        Assume towards a contradiction that \(S\subseteq A\) has the property that
        both \(S\) and \(A\setminus S\) are \(F\)-positive.
        Let \(U\) and \(W\) be \(\epsilon\)-complete ultrafilters
        extending \(F\restriction S\) and \(F\restriction (A\setminus S)\) respectively.
        Fix \(B\) and \(C\) in \(\mathcal P\) such that
        \(B\in U\) and \(C\in W\). Since \(S\subseteq A\) and \(S\in U\),
        \(S\cap B\neq \emptyset\), and it follows that \(A = B\)
        since the sets in \( \mathcal P\) are disjoint. Similarly \(A = C\).
        Since each \(D\in \mathcal P\) belongs to exactly
        one ultrafilter in \(\mathcal T\), \(U = W\),
        which contradicts that \(A\in U\) and \(A\setminus S\in W\).
    \end{proof}
\end{thm}

\begin{thm}\label{thm:measurables}
    If there is a rank Berkeley cardinal,
    then there is a closed unbounded class of cardinals \(\kappa\)
    such that either \(\kappa\) or \(\kappa^+\) is measurable.
    \begin{proof}
        Let \(\lambda\) be the least rank Berkeley cardinal and
        let \(\kappa \geq \lambda\) be an almost supercompact cardinal.
        Let \(F\) be the \(\omega\)-closed unbounded filter on \(\kappa^+\).
        By \cref{thm:completeness}, \(F\) is \(\kappa\)-complete. 
        If \(\alpha > \kappa\) and
        \(j : V_\alpha\to V_\alpha\) is an elementary embedding
        with \(j(\kappa) = \kappa\) and \(\crit(j) < \lambda\),
        then \(j[\kappa^+]\) is \(\omega\)-closed unbounded
        and hence \(j[\kappa^+]\in F\). Therefore
        \(F\) is in the normal fine ultrafilter \(\mathcal U\)
        on \(P(\kappa^+)\) derived from \(j\)
        using \(j[\kappa^+]\). Since \(\kappa\) is almost
        supercompact, \cref{lma:almost_super_crit} and \cref{prp:normal_ultrafilters}
        yield a \(\kappa_\mathcal U^+\)-complete fine
        ultrafilter \(\mathcal W\)
        as in the hypothesis of \cref{thm:saturation},
        and so \(F\) has an atom, which means that
        \(F\) extends to a \(\kappa\)-complete ultrafilter.
        Let \(\delta\) be the completeness
        of \(F\). Then \(\delta\) is a measurable cardinal
        and \(\kappa \leq \delta \leq \kappa^+\). This proves the theorem.
    \end{proof}
\end{thm}
\section{Filters}\label{section:filter_extension}
\subsection{The Ketonen order}\label{KOSection}
The \textit{Ketonen order} is a wellfounded partial order of
countably complete filters introduced in the author's thesis \cite{UA}. 
The wellfoundedness
of the Ketonen order cannot be proved
without appeal to the Axiom of Dependent Choice (DC),
and so to avoid using DC in our applications below,
we work instead with the Ketonen order on \textit{wellfounded} filters,
whose wellfoundedness can be proved in ZF alone.

If \(F\) is a filter on \(X\), \(A\in F\), 
and \(\langle G_x : x\in A\rangle\) is a sequence of filters on \(Y\),
then \(F\ulim_{x\in A} G_x\) denotes the filter on \(Y\)
consisting of all \(B\subseteq Y\)
such that for \(F\)-almost all \(x\in A\), \(B\in G_x\).
Also \(F\usum_{x\in A}G_x\) denotes the filter
on \(X\times Y\) consisting of all \(B\subseteq X\times Y\)
such that for \(F\)-almost all \(x\in A\), for \(G_x\)-almost all \(y\in Y\),
\((x,y)\in B\). (There is really no extra generality in
allowing the index set \(A\) to differ from the underlying set \(X\),
but it is occasionally useful for notational purposes.)

\begin{lma}\label{lma:wellfounded_sum}
    Suppose \(F\) is a wellfounded filter on \(X\), 
    and \(\langle G_x : x\in X\rangle\) is a 
    sequence of wellfounded filters on \(Y\).
    Then \(F\usum_{x\in A}G_x\) and \(F\ulim_{x\in A} G_x\) are wellfounded.
    \begin{proof}
        Notice that for any sequence of structures 
        \(\langle \mathcal M_{x,y}:(x,y)\in X\times Y\rangle\),
        there is a homomorphism from
        \(\left(\prod_{(x,y)\in X\times Y} \mathcal M_{x,y}\right) /\left({\textstyle F\usum_{x\in X}G_x}\right)\)
        into \(\prod_{x\in X}\left(\prod_{y\in Y}\mathcal M_{x,y}/G_x\right)/F\).
        Clearly if \(\mathcal M_{x,y}\) is a wellfounded structure,
        then \(\prod_{x\in X}\left(\prod_{y\in Y}\mathcal M_{x,y}/G_x\right)/F\)
        is wellfounded, and hence so is 
        \(\prod_{(x,y)\in X\times Y} \mathcal M_{x,y}/\left(F\usum_{x\in X}G_x\right)\).
        This proves that \(F\usum_{x\in X}G_x\) is wellfounded.
        But \(F\ulim_{x\in X} G_x\) is the pushforward of \(F\usum_{x\in X}G_x\)
        by the projection to the second coordinate, so it follows that 
        \(F\ulim_{x\in X} G_x\) is wellfounded as well. More concretely,
        for any sequence of structures \(\langle \mathcal M_{y}\rangle_{y\in Y}\),
        there is an embedding from \(\prod_{y\in Y} \mathcal M_y /F\ulim_{x\in X} G_x\)
        into \(\prod_{(x,y)\in X\times Y} \mathcal M_y/F\usum_{x\in X} G_x\)
        induced by the map
        sending \(f\in \prod_{y\in Y} \mathcal M_y\) to \(\tilde f\in \prod_{(x,y)\in X\times Y} \mathcal M_y\)
        where \(\tilde f(x,y) = f(y)\).
    \end{proof}
\end{lma}

For each ordinal \(\delta\), let
\(\mathscr F(\delta)\) denote the set of wellfounded filters on \(\delta\),
and let \(\mathscr F(\delta,\alpha)\) denote the set
of \(F\in \mathscr F(\delta)\) such that \(\alpha\in F\).
For \(F,G\in \mathscr F(\delta)\), define the Ketonen order by setting 
\(F\ke G\) if there is a sequence 
\(\langle F_\alpha :\alpha < \delta\rangle\in \prod_{\alpha < \delta} \mathscr F(\delta,\alpha)\)
such that \(F \subseteq G\ulim_{\alpha < \delta} F_\alpha\).
Similarly, set \(F\keq G\) if there is a sequence 
\(\langle F_\alpha :\alpha < \delta\rangle\in \prod_{\alpha < \delta} \mathscr F(\delta,\alpha+1)\)
\(F \subseteq G\ulim_{\alpha < \delta} F_\alpha\).
\begin{thm}\label{thm:ketonen_wf}
    For every ordinal \(\delta\), Ketonen order on \(\mathscr F(\delta)\) is wellfounded.
    \begin{proof}
        Assume by induction that the Ketonen order on \(\mathscr F(\alpha)\) is wellfounded
        for all \(\alpha < \delta\), and we will show that the Ketonen order on 
        \(\mathscr F(\delta)\) is wellfounded. Note that the Ketonen order on
        \(\mathscr F(\delta,\alpha)\) is isomorphic to the Ketonen order on \(\mathscr F(\alpha)\),
        and so \((\mathscr F(\delta,\alpha),\ke)\) is wellfounded.

        Fix a wellfounded filter \(G\in \mathscr F(\delta)\) such that the Ketonen order is illfounded below \(G\).
        For \(\vec F \in \prod_{\alpha < \delta} \mathscr F(\delta,\alpha)\),
        let \(\xi(\vec F)\) be the rank of \([\vec F]_G\) in the structure 
        \(\prod_{\alpha < \delta} (\mathscr F(\delta,\alpha),\ke)/G\). This structure is wellfounded 
        because both the filter \(G\) and
        the structures \((\mathscr F(\delta,\alpha),\ke)\) are 
        wellfounded.

        For each \(F\ke G\), let \(\xi(F)\) be the minimum ordinal
        of the form \(\xi(\vec F)\) for some
        \(\vec F \in \prod_{\alpha < \delta} \mathscr F(\delta,\alpha)\) such that
        \(F \subseteq G\ulim \vec F\). 
        Then if \(E\ke F\ke G\), we claim \(\xi(E) < \xi(F)\).
        To see this, fix \(\vec F\in \prod \mathscr F(\delta)\) such that \(\xi(F) = \xi(\vec F)\)
        and \(F \subseteq G\ulim_{\alpha < \delta} \vec F\),
        and fix \(\vec E \in \prod \mathscr F(\delta,\alpha)\)
        such that \(E = F\ulim \vec E\). Let 
        \(D_\alpha = F_\alpha\ulim \vec E\)
        where \(\langle F_\alpha\rangle_{\alpha < \delta} = \vec F\). 
        By \cref{lma:wellfounded_sum}, \(D_\alpha\in \mathscr F(\delta,\alpha)\).
        Moreover, 
        \(\vec D= \langle D_\alpha\rangle_{\alpha < \delta}\) 
        is strictly below \(\vec F\) in 
        \(\prod_{\alpha < \delta} \mathscr F(\delta,\alpha)/G\)
        and \(E = G\ulim \vec D\).
        Therefore \(\xi(E) \leq \xi(\vec D) < \xi(\vec F) = \xi(F)\), as claimed.
        It follows that the map \(F\mapsto \xi(F)\) ranks the
        Ketonen order restricted to the predecessors of \(G\), which contradicts
        the assumption that the Ketonen order is illfounded below \(G\).
    \end{proof}
\end{thm}

The following lemmas clarify the relationship between \(\keq\) and \(\ke\)
to a certain extent.
\begin{lma}
Suppose \(F,G\in\mathscr F(\delta)\) satisfy
\(F\keq G\). Then either \(F\ke G\) or there is some \(S\in G^+\) such that
\(F\subseteq G\restriction S\). 
\begin{proof}
Let \(F \subseteq G\ulim_{\alpha < \delta}F_\alpha\)
where \(F_\alpha\in \mathscr F(\delta,\alpha+1)\) for \(\alpha < \delta\). Since
\(F\not\ke G\), the set \(S = \{\alpha < \delta : \alpha\notin F_\alpha\}\) is \(G\)-positive.
For all \(\alpha\in S\), \(F_\alpha\subseteq p_\alpha\) where \(p_\alpha\)
is the principal ultrafilter on \(\delta\) concentrated at \(\alpha\), and so
\[F\subseteq G\ulim_{\alpha < \delta}F_\alpha \subseteq
(G\restriction S)\ulim_{\alpha < \delta}F_\alpha \subseteq 
(G\restriction S)_\ulim_{\alpha < \delta}p_\alpha = G\restriction S\qedhere\]
\end{proof}
\end{lma}
\begin{cor}\label{cor:k_antisymmetric}
    If \(U\keq W\) are ultrafilters in \(\mathscr F(\delta)\), either \(U\ke W\) or \(U = W\).\qed
\end{cor}

One can characterize the Ketonen order as a sort of reducibility on
\(\mathscr F(\delta)\), which turns out to be quite useful.
A function \(f : P(\delta)\to P(\delta)\) is \textit{Lipschitz}
if for all \(\alpha < \delta\), if \(A,B\in P(\delta)\)
satisfy \(A\cap \alpha = B\cap \alpha\), then
\(f(A) \cap \alpha = f(B)\cap \alpha\). 
A Lipschitz function \(f\) is \textit{Ketonen} if
for any \(F\in \mathscr F(\delta)\),
\(f^{-1}[F]\in \mathscr F(\delta)\). 

\begin{lma}\label{lma:ketonen_reduction}
    For \(F,G\in\mathscr F(\delta)\), then \(F\keq G\)
    if and only if there is a Ketonen function \(f : P(\delta)\to P(\delta)\) 
    such that \(F\subseteq f^{-1}[G].\)
    \begin{proof}
        Suppose \(F\keq G\). Fix \(\langle F_\nu:\nu < \delta\rangle\in 
        \prod_{\nu < \delta}\mathscr F(\delta,\nu+1)\) such that
        \(F \subseteq G\ulim_{\nu < \delta} F_\nu\).
        Define \(f : P(\delta)\to P(\delta)\) by 
        \(f(A) = \{\nu < \delta : A\in F_\nu\}\).
        Since for \(\nu < \delta\), \(F_\nu\in F(\delta,\nu+1)\),
        if \(A,B\in P(\delta)\) agree up to some \(\alpha < \delta\), then
        \(\{\nu < \alpha: A\in F_\nu\} = \{\nu < \alpha : B\in F_\nu\}\).
        It follows that \(f\) is Lipschitz. Clearly for all
        filters \(H\) on \(\delta\),
        \(f^{-1}[H] = H\ulim_{\alpha < \delta}F_\alpha\).
        It follows from this and \cref{lma:wellfounded_sum} that
        \(f\) is Ketonen, and moreover \(F \subseteq f^{-1}[G]\).

        Conversely, if \(f : P(\delta)\to P(\delta)\) is a Ketonen function
        such that \(F \subseteq f^{-1}[G]\),
        then define \(\langle F_\nu : \nu < \delta\rangle\) by
        \(F_\nu = \{A \subseteq \delta : \nu \in f(A)\}\). Since
        \(F_\nu = f^{-1}[p_\nu]\) where \(p_\nu\) is the principal ultrafilter
        on \(\delta\) concentrated at \(\nu\),
        \(F_\nu \in \mathscr F(\delta)\). 
        Since \(f\) is Lipschitz and
        \(\nu \in f(\delta)\), \(\nu \in f(\nu+1)\),
        and hence \(F_\nu\in \mathscr F(\delta,\nu+1)\). 
        Moreover,
        for any \(A\subseteq \delta\),
        \(f(A) = \{\nu < \delta : A\in F_\nu\}\).
        Hence \(F\subseteq f^{-1}[G] = G\ulim_{\nu < \delta}F_\nu\).
        This shows that \(F\keq G\).
    \end{proof}
\end{lma}

For any filter \(F\in \mathscr F(\delta)\), \(|F|_{\Bbbk}\) denotes the rank
of \(F\) in the Ketonen order on \(\mathscr F(\delta)\).
\begin{lma}\label{lma:fixer}
    Suppose \(U\in \mathscr F(\delta)\) is an
    ultrafilter, \(V_\alpha\preceq_{\Sigma_1} V\),
    and \(j : V_\alpha\to V_\alpha\)
    is an elementary embedding fixing \(\delta\) and \(|U|_\Bbbk\).
    Then \(j(U) = U\).
    \begin{proof}
        For any \(F\in \mathscr F(\delta)\), 
        \(V_\alpha\) correctly computes
        \(|F|_{\Bbbk}\):
        since \(V_\alpha\preceq_{\Sigma_1} V\),
        \(V_\alpha\) correctly computes \(\mathscr F(\delta)\),
        and so \(V_\alpha\) is correct about the structure
        \((\mathscr F(\delta),\ke)\), and the fact that 
        \(V_\alpha\preceq_{\Sigma_1} V\)
        also suffices to conclude that \(V_\alpha\) correctly computes 
        the rank function associated to any wellfounded relation
        in \(V_\alpha\).

        By \cref{lma:ketonen_reduction}, \(U \keq j(U)\): \(f = j\restriction P(\delta)\)
        is a Ketonen function and \(U \subseteq f^{-1}[j(U)]\).
        On the other hand, \(U\not\ke j(U)\) since the Ketonen rank of \(j(U)\)
        is \(j(|U|_\Bbbk) = |U|_\Bbbk\), the same as the Ketonen rank of \(U\).
        Therefore by \cref{cor:k_antisymmetric}, \(j(U) = U\).
    \end{proof}
\end{lma}

\subsection{Interlude: the theory of \(j\)-of-\(j\)-of-\(j\)}
If \(j : V_\alpha\to V_\alpha\) is an elementary embedding
and \(\alpha\) is a limit ordinal, we define
the iterates \(\langle j_n:n < \omega\rangle\) of \(j\) by
recursion, setting \(j_0 = j\), and for each \(n < \omega\),
\(j_{n+1} = j(j_n)\). We include two folklore
observations about iterates which are closely
related to Kunen's proof of the wellfoundedness of iterated ultrapowers
\cite{KunenLU}.

First, for \(\beta < \alpha\), \(j_n(\beta) \geq j_{n+1}(\beta)\) with equality
if and only if \(\beta \in \ran( j_n)\). To see this, first note
that for all
\(\gamma < \alpha\), we have \(j_n(\sup j_n[\gamma]) = \sup j_{n+1}[j_n(\gamma)]\).
Then let \(\gamma < \alpha\) be least such that \(j_n(\gamma) \geq \beta\),
so that \(\sup j_n[\gamma]\leq \beta \leq j_n(\gamma)\). If \(\beta < j_n(\gamma)\), then \(\gamma\)
is a limit ordinal and \(j_{n+1}(\beta) < \sup j_{n+1}[j_n(\gamma)] = j_n(\sup j_n[\gamma]) \leq j_n(\beta)\).
Otherwise \(\beta = j_n(\gamma)\), so \(\beta\in \ran(j_n)\) and \(j_n(\beta) = 
j_n(j_n(\gamma)) = j_{n}(j_n)(j_n(\gamma)) = j_{n+1}(\beta)\).

Second, for any ordinal \(\beta < \alpha\),
there is some \(n < \omega\) such that \(j_n(\beta) = \beta\).
Towards a contradiction, suppose \(\beta < \alpha\) is the least counterexample.
The sequence \(\langle j_n(\beta) : n < \omega\rangle\)
is weakly decreasing by the previous paragraph,
so fix \(n < \omega\) such that \(j_n(\beta) = j_{n+1}(\beta)\). 
Then \(\beta = j_n(\gamma)\) for some \(\gamma \leq \beta\).
If \(\gamma = \beta\), we are done, so assume that \(\gamma < \beta\).
Then by the minimality of \(\beta\), there is some
\(m > n\) such that \(j_{m}(\gamma) = \gamma\). Then 
\(j_{m+1}(\beta) = j_{n}(j_m)(\beta) = j_{n}(j_m)(j_n(\gamma)) = j_{n}(j_{m}(\gamma)) = j_n(\gamma) = \beta\).
\subsection{The semi-linearity of the Ketonen order}\label{section:semilinear}
In the context of the Axiom of Choice, the \textit{Ultrapower Axiom} (UA) roughly states
that any two ultrapowers of the universe have a common internal ultrapower.
The principle has found a number of applications in the theory of supercompact cardinals.
UA is equivalent to the linearity of the Ketonen order on ultrafilters.\footnote{
    This is defined for countably complete
    ultrafilters \(U,W\) on \(\delta\) by setting \(U < W\) if
    \(U = W\ulim_{\alpha < \delta} U_\alpha\) 
    where for \(W\)-almost all \(\alpha < \delta\),
    \(U_\alpha\) is a countably complete ultrafilter on \(\delta\) 
    concentrating on \(\alpha\).
}
Here we will show that in the context of a rank Berkeley cardinal, the Ketonen order is
\textit{almost linear} in the sense that it contains no large antichains. 
This fact will find an application in the proof of the filter extension
property below.

One can actually show the
semilinearity of a finer order than the Ketonen order.
The \textit{embedding order}
is defined by setting \(F \eeq G\) if \(F\subseteq j^{-1}[G]\) for a Ketonen 
\textit{elementary} embedding
\(j : P(\delta)\to P(\delta)\). Here we view \(P(\delta)\) as a transitive structure
endowed with the membership relation.

\begin{thm}\label{thm:antichains}
    Suppose \(\lambda\) is the least
    rank Berkeley cardinal
    and \(A\subseteq \mathscr F(\delta)\) is a set of ultrafilters that
    are incomparable in the embedding order.
    \begin{enumerate}[(1)]
        \item\label{incomp} If \(\kappa\) is almost supercompact and
        \(\cf(\kappa) > \lambda\), then \(\scott(A) < \kappa\).
        \item\label{woincomp} If \(A\) can be wellordered, then \(|A|\leq \lambda\).
    \end{enumerate}
    \begin{proof}
        The key observation is that if \(U\in
        \mathscr F(\delta)\) is an ultrafilter 
        and \(j : V_{\alpha}\to V_{\alpha}\)
        is an elementary embedding where \(\alpha > \delta\)
        and \(V_\alpha\preceq_{\Sigma_1} V\), 
        then for some \(n < \omega\),
        \(j_n(U) = U\). To see this, find \(n < \omega\) such that
        \(j_n(\xi) = \xi\) where \(\xi\) is the rank of \(U\) in 
        the Ketonen order on \(\mathscr F(\delta)\), and apply \cref{lma:fixer}.

        Fix \(\alpha\) larger than the rank of \(A\) such 
        that \(V_\alpha\preceq_{\Sigma_1} V\),
        and suppose \(j : V_\alpha\to V_\alpha\) is an elementary embedding
        with \(\kappa_\omega(j) = \lambda\).
        For each \(n < \omega\), let \(A_n\) be the set of ultrafilters
        in \(A\) fixed by \(j_n\). Then \(A = \bigcup_{n < \omega} A_n\).
        
        We claim that \(j_n(A_n) = j_n[A_n]\).
        (Of course, \(j_n[A_n] = A_n\).) 
        Suppose not, and fix \(U\in j_n(A_n)\setminus j_n[A_n]\).
        Let \[B = j_{n+1}(j_n(A_n))
        = j_n(j_n(A_n))\] 
        Then \(j_n(U)\) and \(j_{n+1}(U)\) belong to \(B\),
        and \(B\) is a set of Ketonen incomparable ultrafilters.
        Since every ultrafilter in \(A_n\) is fixed by \(j_n\),
        every ultrafilter in \(j_n(A_n)\) is fixed by \(j_{n+1}\),
        and it follows that \(j_{n+1}(U) = U\).
        On the other hand, clearly
        \(U\eeq j_n(U)\). 
        Also \(U\neq j_n(U)\) since
        \(U\notin j_n[A_n]\). This contradicts that \(B\)
        is a set of incomparable elements of the embedding
        order on \(\mathscr F(\delta)\).
        
        We now prove \ref{incomp}.
        By \cref{lma:j_closed}, since \(j_n(A_n) = j_n[A_n]\), 
        there is some \(\beta_n < \kappa\) such that
        \(A_n\) injects into \(V_{\beta_n}\).
        Applying the wellordered collection lemma,
        there is some \(\gamma < \kappa\) and 
        a set \(\{f_x : x\in V_\gamma\}\) such that for each \(n\), there is
        some \(x\in V_\gamma\) such that
        \(f_x : A_n\to V_{\beta_n}\) is an injection.
        Define \(g_n : A_n\to {}^{V_\gamma} V_{\beta_n}\)
        by \(g_n(a)(x) = f_x(a)\). Then \(g_n\) is an injection
        from \(A_n\) into \(V_{\rho_n}\) where
        \(\rho_n = \gamma \cdot \beta_n + 1\). 
        Now setting \(g(a) = \langle g_n(a) : n < \omega\rangle\),
        we obtain an injection from \(A\) to \(V_\beta\)
        where \(\beta = (\sup_{n < \omega} \rho_n) + 1\) is less than \(\kappa\).

        Assuming \(A\) can be wellordered, one can conclude from the
        fact that \(j_n(A_n) = j_n[A_n]\) that \(|A_n| < \kappa_n(j)\),
        and hence \(|A| \leq \lambda\), proving \ref{woincomp}.
    \end{proof}
\end{thm}
It is unclear whether \cref{thm:antichains} \ref{incomp} can be improved to 
\(|A| < \lambda\).

The semi-linearity property of the 
Ketonen order that will actually be applied in the proof of the filter
extension property is more technical
but a bit easier to show.
\begin{prp}\label{prp:rank_closed}
    Suppose \(\lambda\) is a rank Berkeley cardinal, \(\delta\) and \(\xi\)
    are ordinals,
    and \(\mathscr U\) is the set of all ultrafilters
    of rank \(\xi\) in the Ketonen order on \(\mathscr F(\delta)\). 
    Then \(\lambda\) is \(\mathscr U\)-closed rank Berkeley. 
    \begin{proof}
        Suppose \(\alpha > \delta\) is an ordinal such that \(V_\alpha\preceq_{\Sigma_1} V\)
        and \(j : V_\alpha\to V_\alpha\) is an elementary embedding that fixes \(\delta\) and
        \(\xi\). Then \(j(\mathscr U) = \mathscr U\) and for any 
        ultrafilter \(W\in \mathscr U\), \(j(W) = W\) by \cref{lma:fixer}.
        It follows that \(j[\mathscr U] = \mathscr U = j(\mathscr U)\).
        This easily implies that \(\lambda\) is \(\mathscr U\)-closed rank Berkeley.
    \end{proof}
\end{prp}

\subsection{Atoms of the closed unbounded filter}\label{section:closed unboundeds}
In this section, we use the theory of the Ketonen order to give a deeper
analysis of the structure of the closed unbounded filter.
In particular, we study the atoms of the closed unbounded filter
and their relationship with stationary reflection.
\begin{lma}
    Suppose \(F\) and \(G\) are extensions of the closed unbounded filter
    on a regular cardinal \(\delta\). If \(F\eeq G\),
    either \(F\subseteq G\) or for all \(S\in F\), 
    for a \(G\)-positive set of \(\alpha < \delta\), \(S\cap \alpha\) is stationary.
    \begin{proof}
        Let \(j : P(\delta)\to P(\delta)\) be an elementary embedding
        such that \(F \subseteq j^{-1}[G]\).
        
        First consider the case that 
        the set \(E = \{\alpha < \delta : j(\alpha) = \alpha\}\)
        belongs to \(G\). We claim that \(F\subseteq G\). To see this 
        fix \(A\in F\). Since \(F\subseteq j^{-1}[G]\),
        \(j(A)\in G\), and so \(j(A)\cap E\in G\), which implies that 
        \(A\in G\) because \(j(A)\cap E = A\cap E \subseteq A\).

        Assume instead that \(\delta\setminus E\) is \(G\)-positive.
        Fix \(S\in F\).
        Let \(C\) be the set of closure points of \(j\).
        Then \(C\in G\) since \(G\) extends the closed unbounded filter.
        Fix \(\alpha\in j(S)\cap C\setminus E\), 
        and we will show that \(S\cap \alpha\) is stationary.
        Fix a closed unbounded set \(B\subseteq \alpha\),
        and we will show that \(S\cap B\) is nonempty. Then \(j(B)\)
        is closed unbounded in \(j(\alpha)\)
        and \(j[\alpha]\subseteq j(B)\).
        Since \(\alpha = \sup j[\alpha]\) is a limit point of \(j(B)\)
        and \(\alpha < j(\alpha)\) since \(\alpha\notin E\),
        \(\alpha\) must belong to the closed unbounded set \(j(B)\). 
        Since \(\alpha\in j(S)\cap j(B)\),
        by elementarity, \(S\cap B\) is nonempty.

        We have shown that for all \(\alpha\in j(S)\cap C\setminus E\),
        \(S\cap \alpha\) is stationary. Since \(j(S)\cap C\in G\)
        and \(\delta\setminus E\) is \(G\)-positive,
        \(j(S)\cap C\setminus E\) is \(G\)-positive, and so we are done.
    \end{proof}
\end{lma}

If \(S,T\subseteq \delta\)
are stationary sets, \(S\subseteq \delta\) \textit{reflects stationarily in \(T\)} 
if there is a stationary set
of \(\alpha\in T\) such that \(S\cap \alpha\) is stationary in \(\alpha\);
\(S\) \textit{reflects fully in \(T\)} if for all but a nonstationary
set of \(\alpha\in T\), \(S\cap \alpha\) is stationary.
No stationary set can reflect fully in itself. 
A stationary set is \textit{thin} if it does not even reflect stationarily in itself.
\begin{cor}\label{cor:club_incomparable}
    Suppose \(S\) is a thin stationary
    subset of a regular cardinal \(\delta\).
    Let \(D\) be the closed unbounded filter on \(\delta\) restricted to \(S\).
    Then if \(F\) and \(G\) extend \(D\),
    \(F\eeq G\) if and only if \(F\subseteq G\).\qed
\end{cor}
Any atom of the closed unbounded filter is thin, since otherwise
it would reflect fully in itself.
Under choiceless large cardinal axioms, a partial converse holds:
\begin{thm}\label{thm:atoms}
    Suppose \(\lambda\) is rank Berkeley.
    Then for all sufficiently large regular cardinals \(\delta\),
    any thin stationary subset of \(\delta\)
    is the disjoint union of at most \(\lambda\)
    atoms of the closed unbounded filter.
    \begin{proof}
        Let \(\kappa\geq\lambda\) be almost supercompact
        and suppose \(\delta \geq \kappa\) is a 
        regular cardinal large enough that the closed
        unbounded filter \(\mathcal C\) on \(\delta\) is \(\aleph(V_{\kappa+1})\)-complete.\footnote{
            One can show that the least such
            \(\kappa\) is singular. By \cref{thm:measurables},
            \(\kappa^+\) is measurable,
            which implies \(\aleph(V_{\kappa+1}) = \kappa^+\),
            and so by \cref{thm:completeness},
            the closed unbounded filter on \textit{any} regular
            \(\delta\geq \kappa\) is \(\aleph(V_{\kappa+1})\)-complete.
        }

        Since \(\mathcal C\) is
        \(\kappa\)-complete by \cref{thm:completeness},
        the proof of \cref{lma:complete_wf} yields a wellorderable
        family of ultrafilters 
        \(\mathcal D\subseteq \mathscr F(\delta)\) such that
        for any stationary set \(T\subseteq \delta\),
        there is an ultrafilter \(U\in \mathcal D\) that
        extends \(\mathcal C\restriction T\).

        Let \(\mathcal A\) be the set of all \(U\in \mathcal D\)
        extending \(\mathcal C\restriction S\). 
        By \cref{cor:club_incomparable},
        \(\mathcal A\) is a set of incomparable ultrafilters in the embedding order
        on \(\mathscr F(\delta)\), and so by \cref{thm:antichains},
        \(|\mathcal A|\leq \lambda\). 

        The remainder of the proof is similar to
        \cref{lma:saturation}.
        By the wellordered collection lemma,
        there is a set \(\sigma\) such that \(\dscott(\sigma) \leq \kappa\) 
        and for any \(U_0,U_1\in \mathcal A\), there are disjoint \(T_0,T_1\in \sigma\)
        such that \(T_0\in U_0\) and \(T_1\in U_1\).
        For each \(\alpha < \delta\), let \(D_\alpha = \{T\in \sigma : \alpha\in T\}\)
        and for \(\alpha,\beta \in S\), set \(\alpha\sim \beta\)
        if \(D_\alpha = D_\beta\).
        Then the equivalence classes of \(\sim\) are in one-to-one correspondence
        with the wellorderable family 
        \(\{D_\alpha : \alpha \in S\}\subseteq P(\sigma)\).
        It follows that \(|S/{\sim}| < \aleph(V_{\kappa+1})\). 

        Let \(\mathcal T\subseteq P(S)\) be the set of stationary 
        \(\sim\)-equivalence classes.
        Since \(\bigcup (S/{\sim}) = S\) and 
        \(\mathcal C\restriction S\) is \(\aleph(V_{\kappa+1})\)-complete,
        \(\bigcup \mathcal T\in \mathcal C\restriction S\). 
        In other words, \(S\setminus \bigcup \mathcal T\) is nonstationary.
        To finish, we will show that if \(T\in \mathcal T\), then \(T\) is an atom of
        \(\mathcal C\). 
        
        First, note that \(T\)
        belongs to at most one ultrafilter in \(\mathcal A\):
        if \(T\in U\in \mathcal A\), then for any \(R\in U\cap \sigma\),
        \(T\subseteq R\): fix any \(\alpha\in T\cap R\), and
        note that if \(\beta\in T\), then since \(\alpha\sim \beta\),
        \(\beta\in R\). Therefore by our choice of \(\sigma\),
        \(T\) is \(W\)-null for all \(W\in \mathcal A\) other than \(U\).

        Finally, assume towards a contradiction that \(T\) is not
        an atom, and so \(T\) is the 
        disjoint union of two stationary sets \(B_0\) and \(B_1\). Then there exist
        \(U_0,U_1\in \mathcal A\) such that \(B_0\in U_0\) and \(B_1\in U_1\).
        But then \(T\in U_0\cap U_1\) and \(U_0\neq U_1\), 
        which contradicts our observation that
        \(T\) belongs to at most one \(U\in \mathcal A\). 
    \end{proof}
\end{thm}
\begin{thm}
    Suppose \(\lambda\) is rank Berkeley.
    Then for all sufficiently large regular cardinals \(\delta\), 
    the closed unbounded filter on \(\delta\) is atomic.
    \begin{proof}
        If \(T\subseteq \delta\) is stationary,
        let \(S\) be the set of \(\alpha\in T\) such that
        \(T\cap \alpha\) is nonstationary. Since \(T\) does not
        reflect fully in itself, \(S\) is stationary, and clearly
        \(S\) is thin. By \cref{thm:atoms} \(S\) contains an atom,
        and since \(S\subseteq T\), so does \(T\).
    \end{proof}
\end{thm}
The \textit{reflection order} is defined on stationary sets
\(S,T\subseteq \delta\) by setting \(S < T\) if \(S\) reflects fully in \(T\). 
Two stationary sets are \textit{reflection incomparable}
if they are incomparable in the reflection order.
The following theorem states that restricted to atoms,
the reflection order is almost linear:
\begin{thm}\label{thm:reflection_order}
    Suppose \(\lambda\) is rank Berkeley. 
    Then for all sufficiently large regular cardinals
    \(\delta\), any set of reflection incomparable
    atoms of the closed unbounded filter on \(\delta\) has cardinality
    less than or equal to \(\lambda\).
    \begin{proof}
        Suppose \(\kappa \geq \lambda\) is almost
        supercompact and \(\delta \geq \kappa\) is a regular cardinal.
        Let \(\mathcal C\) denote the closed unbounded filter on \(\delta\).
        Suppose \(\mathcal S\) is a family of atoms
        of \(\mathcal C\) that are incomparable in the
        reflection order. 
        Let \(\mathcal A = \{\mathcal C\restriction T : T\in \mathcal S\}\).
        By \cref{cor:club_incomparable}, 
        \(\mathcal A\) is a set of ultrafilters that are incomparable
        in the embedding order. By \cref{thm:completeness}, every \(U\in \mathcal A\)
        is \(\kappa\)-complete, and so by \cref{lma:complete_wf}, 
        \(\mathcal A\) is wellorderable.
        It follows that \(|\mathcal A| \leq \lambda\).
    \end{proof}
\end{thm}
The proof shows the stronger statement 
that every family of reflection incomparable \textit{equivalence classes} of
atoms of the closed unbounded filter modulo the nonstationary ideal 
has cardinality less than or equal to \(\lambda\).
\subsection{The filter extension property}
Suppose there is a rank Berkeley cardinal and \(\kappa\) is almost extendible.
Can every \(\kappa\)-complete filter on an ordinal
be extended to a \(\kappa\)-complete ultrafilter? 
\cref{lma:complete_wf}
shows that any \(\kappa\)-complete filter on an ordinal can be extended
to an ultrafilter, but controlling its completeness 
requires much more effort. The problem is that the ultrafilters
used in \cref{lma:complete_wf} are derived from extendibility embeddings;
such an ultrafilter is a priori only as complete as the critical point of the embedding
from which it is derived, and under our hypotheses, there may be no extendibility embeddings
with critical point larger than the least rank Berkeley cardinal.
The ultrafilters we will use instead come from \cref{thm:saturation}. 
\begin{thm}\label{thm:filter_extension}
    Assume there is a rank Berkeley cardinal.
    Then for a closed unbounded class of cardinals \(\kappa\), every \(\kappa\)-complete
    filter on an ordinal extends to a \(\kappa\)-complete ultrafilter.
    \begin{proof}
        Let \(\lambda\) be the least rank Berkeley cardinal
        and let \(\Gamma\) be the class of all \(X\) such that
        \(\lambda\) is \(X\)-closed rank Berkeley.
        Let \(\kappa\) be a cardinal that 
        is \(X\)-closed almost extendible for 
        all \(X\in \Gamma\); by \cref{thm:closed_extendible}, there is a closed unbounded class of
        such cardinals. 
        We will show that every \(\kappa\)-complete filter on an ordinal 
        extends to a \(\kappa\)-complete ultrafilter.

        Fix an ordinal \(\delta\) and assume towards a contradiction
        that there is a \(\kappa\)-complete filter on \(\delta\) that
        does not extend to a \(\kappa\)-complete ultrafilter.
        By \cref{lma:complete_wf}, every \(\kappa\)-complete filter
        on \(\delta\) extends to a wellfounded
        ultrafilter and hence is itself wellfounded.
        In particular, every \(\kappa\)-complete filter on \(\delta\) 
        belongs to \(\mathscr F(\delta)\). Let \(\mathscr S\) be the set
        of \(\kappa\)-complete filters on \(\delta\) that
        do not extend to \(\kappa\)-complete ultrafilters. Applying \cref{thm:ketonen_wf},
        let \(F\) be a Ketonen minimal element of \(\mathscr S\).

        Fix \(\alpha > \delta\) be such that \(V_\alpha\preceq_{\Sigma_1} V\), and
        let \(\mathcal E\) be the set of elementary embeddings
        \(j : V_\alpha\to V_\alpha\) such that \(F\in \ran(j)\).
        (The next paragraph shows that \(\mathcal E\) is nonempty.)
        Let \(\kappa' = \kappa\) if \(\kappa\) is regular and \(\kappa' = \kappa+1\) if \(\kappa\) is singular.
        Let \[\mathcal B = \left\{\bigcap_{j\in \sigma} j[\delta] : 
        \sigma\subseteq\mathcal E,\, \dscott(\sigma) < \kappa'\right\}\]
        Let \(G\) be the filter on \(\delta\) generated by \(\mathcal B\).
        The filter \(G\) is \(\kappa'\)-complete by the wellordered collection lemma
        (\cref{thm:wo_collection}).
        The main claim of the proof is that \(G\cup F\) generates a proper filter,
        or in other words, that every set in \(\mathcal B\) is \(F\)-positive.

        Before proving the claim, let us show how to use it to complete the proof of the theorem.
        Let \(i : V_{\alpha+\omega}\to V_{\alpha+\omega}\) be an elementary embedding
        fixing \(\delta\) such that \(\crit(i) < \lambda\), and note that \(i[\delta]\in i(G)\):
        \(i \restriction V_\alpha\) is an elementary embedding with \(i(F)\in \ran(i)\),
        which means \(i\in i(\mathcal E)\), so \(i[\delta]\in \mathcal B\subseteq G\).
        Therefore \(G\in \mathcal U\) where \(\mathcal U\) is the normal fine ultrafilter
        on \(P(\delta)\) derived from \(i\) using \(i[\delta]\).
        Since \(G\) is \(\kappa\)-complete, \cref{thm:saturation} implies
        that there is a partition of \(\delta\) into fewer than \(\lambda\)-many 
        \(G\)-positive sets
        \(\langle S_\nu : \nu < \gamma\rangle\) such that 
        \(G\restriction S_\nu\) is an ultrafilter for all \(\nu < \gamma\).
        (The hypotheses
        of \cref{thm:saturation}, namely the existence of \(\mathcal W\), hold 
        for sufficiently large \(\delta\) by the remarks following \cref{lma:complete_wf}.)

        Assume towards a contradiction that for all \(\nu < \gamma\),
        \(F\nsubseteq G\restriction S_\nu\). Therefore there is a
        set in \(F\) that is not in \(G\restriction S_\nu\), and so since \(G\restriction S_\nu\)
        is an ultrafilter, there is a \(G\)-large set whose intersection with
        \(S_\nu\) is in the dual ideal \(F^*\). By the wellordered collection lemma,
        there is a set \(\upsilon\subseteq P(\delta)\) such that \(\dscott(\upsilon) < \kappa'\) and
        for each \(\nu < \gamma\), there is a set \(A\in G\cap \upsilon\) such that
        \(A\cap S_\nu\in F^*\). Let \(B = \bigcap \upsilon\). Then \(B\in G\)
        by definition, but \(B\cap S_\nu\in F^*\) for all \(\nu < \gamma\).
        Since \(\langle S_\nu : \nu < \gamma\rangle\) is a partition of
        \(\delta\), \(B = \bigcup_{\nu < \gamma} B\cap S_\nu\).
        Since \(F\) is \(\kappa\)-complete, it follows that
        \(B \in F^*\). Now \(B\in G\) and \(B\) is not \(F\)-positive, contrary to the claim.
        This contradiction establishes that for some \(\nu < \gamma\),
        \(F\subseteq G\restriction S_\nu\), which shows that \(F\) extends
        to a \(\kappa\)-complete ultrafilter.

        We now proceed to the proof of the claim.
        Assume towards a contradiction that for some \(\sigma\subseteq \mathcal E\) with 
        \(\dscott(\sigma) < \kappa'\),
        \(\bigcap_{j\in \sigma}j[\delta]\) is \(F\)-null.
        In other words, 
        \[S = \bigcup_{j\in \sigma}(\delta\setminus j[\delta]) \in F\]
        For each \(\xi < \delta\) and \(j\in \sigma\), let \(D_\xi(j)\)
        be the ultrafilter on \(\delta\) derived from \(j\) using \(\xi\).
        For all \(\xi\in S\), there is some \(j\in \sigma\) such that
        \(\xi\in D_\xi(j)\) since \(\xi < j(\xi)\).
        Also \(D_\xi(j)\in \mathscr F_\xi(\delta)\) since \(j\)
        is an embedding of \(V_\alpha\) and \(V_\alpha\preceq_{\Sigma_1} V\).

        For \(\xi < \delta\), let \[D_\xi = \bigcap \{D_\xi(j) : j\in \sigma\}\]
        and note that if \(\xi \in S\), then \(\xi \in (D_\xi)^+\) since for some \(j\in \sigma\),
        \(\xi\in D_\xi(j)\).
        For each \(j \in \sigma\), let \(F_j = j^{-1}(F)\).
        Since \(F\in \ran(j)\), \(F_j\) is a \(\kappa\)-complete filter
        that does not extend to a \(\kappa\)-complete ultrafilter.
        We claim that 
        \[\bigcap_{j\in \sigma} F_j \subseteq F\ulim_{\xi < \delta}D_\xi\]
        In fact, equality holds, but we do not need this.

        Suppose \(A\in \bigcap_{j\in \sigma} F_j\).
        Then for all \(j\in \sigma\), 
        \(j(A)\in F\), or in other words, \[\{\xi < \delta: A\in D_\xi(j)\}\in F\]
        Since \(F\) is \(\kappa'\)-complete and \(\dscott(\sigma) < \kappa'\),
        by \cref{lma:ordinal_complete},
        \[\left\{\xi < \delta : A\in \bigcap_{j\in \sigma} D_\xi(j)\right\} \in F\]
        and so 
        \[\{\xi < \delta : A\in D_\xi\} \in F\]
        which implies that \(A\in F\ulim_{\xi < \delta}D_\xi\).

        It follows that \(\bigcap_{j\in \sigma} F_j\ke F\).
        To see this, note that \(S\in F\), and for all \(\xi\in S\),
        \(\xi\in (D_\xi)^+\). Let \(E_\xi = D_\xi \restriction \xi\).
        Then \(E_\xi \in \mathscr F(\delta,\xi)\):
        this follows from the fact that for any \(j\in \sigma\) such that \(j(\xi) > \xi\),
        \(D_\xi(j)\in \mathscr F(\delta,\xi)\) and \(E_\xi\subseteq D_\xi(j)\).
        As a consequence, \(\bigcap_{j\in \sigma} F_j\subseteq F\ulim_{\xi\in S} E_\xi\),
        and so \(\bigcap_{j\in \sigma} F_j\ke F\).

        The intersection \(\bigcap_{j\in \sigma} F_j\) of the \(\kappa\)-complete
        filters \(F_j\) is \(\kappa\)-complete, so
        by the minimality of \(F\), \(\bigcap_{j\in \sigma} F_j\) 
        extends to a \(\kappa\)-complete ultrafilter \(W\).
        We will show that \(F_j\subseteq W\)
        for some \(j\in \sigma\), contradicting that \(F_j\) does not
        extend to a \(\kappa\)-complete ultrafilter.
        This conclusion would be obvious under the Axiom of Choice: 
        if \(F_j\nsubseteq W\) for all \(j\in \sigma\), then for each
        \(j\in \sigma\), choose \(A_j\in F_j\setminus W\), and let 
        \(A = \bigcup_{j\in \sigma} A_j\); then \(A\in \bigcap_{j \in \sigma} F_j\)
        and \(A\notin W\), which contradicts that \(\bigcap_{j\in \sigma} F_j\subseteq W\).
        Since we are working in ZF, a different argument is required.

        \begin{obs}\label{clm}
        If \(\gamma < \kappa'\) and \(\langle H_\alpha\rangle_{\alpha < \gamma}\)
        is a sequence of filters with \(\bigcap_{\alpha < \gamma} H_\alpha\subseteq W\),
        then \(H_\alpha\subseteq W\) for some \(\alpha < \gamma\).
        \end{obs}
        \begin{proof}
            The wellordered collection yields a set \(\sigma\subseteq P(\delta)\)
        such that \(\dscott(\sigma) < \kappa'\) and \(\sigma\cap (H_\alpha\setminus W)\neq \emptyset\)
        for all \(\alpha < \gamma\), and so \(\bigcup (\sigma\setminus W)\)
        belongs to \(\bigcap_{\alpha < \gamma} H_\alpha\) but not \(W\).
        \end{proof}
        We will use this
        observation repeatedly to replace \(\bigcap_{j\in\sigma} F_j\) with a more manageable
        intersection that is still contained in \(W\).

        Our first step is to reduce to the case that \(\dscott(\sigma) < \kappa\). 
        If \(\kappa\) is 
        regular, this is true by definition, so assume instead that \(\kappa\) is singular.
        Since \(\dscott(\sigma)\leq \kappa\),
        one can write \(\sigma = \bigcup_{\alpha < \iota} S_\alpha\) where \(\iota = \cf(\kappa)\)
        and \(\dscott(S_\alpha) < \kappa\) for all \(\alpha < \iota\). 
        Then \(\bigcap_{\alpha <\iota}\bigcap_{j\in S_\alpha}F_j = \bigcap_{j\in \sigma}F_j\subseteq W\),
        and so by \cref{clm}, for some \(\alpha < \iota\),
        \(\bigcap_{j\in S_\alpha}F_j\subseteq W\). By replacing \(\sigma\)
        with \(S_\alpha\), we may assume that \(\dscott(\sigma) < \kappa\).

        For \(\xi < \delta\) and \(j,k\in \sigma\), set \(j\preceq_\xi k\) if
        \(|D_\xi(j)|_\Bbbk \leq |D_\xi(k)|_{\Bbbk}\)
        and \(j\simeq_\xi k\) if \(D_\xi(j) = D_\xi(k)\).
        Let \(Z\) be the set of pairs \(({\preceq},{\simeq})\) where \(\preceq\) is a prewellorder of \(\sigma\)
        and \(\simeq\) is an equivalence relation on \(\sigma\), and
        let \[A({\preceq},{\simeq}) = \{\xi < \delta : ({\preceq}_\xi,{\simeq}_\xi) = {({\preceq},{\simeq})}\}\]
        for any \(({\preceq},{\simeq})\in Z\).
        Then \(\mathcal P = \{A({\preceq},{\simeq}) : ({\preceq},{\simeq})\in Z\}\) is a partition
        of \(\delta\), and therefore it is wellorderable: set \(A < B\) if \(\min(A) < \min(B)\).
        Since \(\dscott(\sigma) < \kappa\) and \(\kappa\) is a limit ordinal, 
        \(\dscott(Z) < \kappa\) as well, and so \(|\mathcal P| < \kappa\).
        (This is the main reason we needed to ensure that \(\dscott(\sigma) < \kappa\).)
        Let \(\mathcal Q = \mathcal P\cap F^+\), and note that
        \(\bigcup \mathcal Q\in F\) since 
        its complement is \(\bigcup (\mathcal P\cap F^*)\)
        which belongs to \(F^*\) by \(\kappa\)-completeness.
        
        Note that in general if
        \(F\) is a \(\kappa\)-complete filter on \(X\) and
        \(\mathcal Q\) is a partition
        of an \(F\)-large set into fewer than \(\kappa\)-many \(F\)-positive sets,
        then for any sequence of filters \(\langle D_x : x\in X\rangle\),
        \[F\ulim_{x\in X}D_x = \bigcap_{A\in \mathcal Q} ((F\restriction A) \ulim_{x\in X} D_x)\]
        Our original argument made use of the wellordered collection lemma, 
        but the anonymous referee suggested the following simpler argument that works in ZF alone.
        For the nontrivial direction, 
        suppose \(S\in \bigcap_{A\in \mathcal Q} ((F\restriction A) \ulim_{x\in X} D_x)\).
        In other words, for each \(A\in \mathcal Q\),
        \(A' = \{x \in A : S\in D_x\}\in F\restriction A\).
        Let \(A'' = (X\setminus A)\cup A'\), so that
        \(A''\in F\). Then \(\bigcap_{A\in Q}A'' \in F\) since \(|Q| < \kappa\)
        and \(F\) is \(\kappa\)-complete, but \(\bigcap_{A\in Q}A'' = \bigcup_{A\in Q} A' = \{x\in X : S\in D_x\}\).
        Since \(\{x\in X : S\in D_x\}\in F\), \(S\in F\ulim_{x\in X}D_x\), as desired.

        Given this, we now have:
        \begin{align*}
            \bigcap_{j\in \sigma} F_j &= \bigcap_{j\in \sigma} j^{-1}[F]\\ 
                                      &= \bigcap_{j\in \sigma} F\ulim_{\xi < \delta} D_\xi(j)\\
                                      &= \bigcap_{j\in \sigma} \bigcap_{A\in \mathcal Q} ((F\restriction A)\ulim_{\xi < \delta} D_\xi(j))\\
                                      &= \bigcap_{A\in \mathcal Q} \bigcap_{j\in \sigma} (F\restriction A)\ulim_{\xi < \delta} D_\xi(j)
        \end{align*}
        We use here that \(F\ulim_{\xi < \delta} D_\xi(j) = j^{-1}[F]\)
        simply because \(j(A) = \{\xi < \delta : A\in D_\xi(j)\}\),
        and so \(A\in j^{-1}[F]\) if and only if \(\{\xi < \delta : A\in D_\xi(j)\}\in F\).
        
        Let \(H_A = \bigcap_{j\in \sigma} (F\restriction A)\ulim_{\xi < \delta} D_\xi(j)\).
        We have proved that \(\bigcap_{A\in \mathcal Q} H_A\subseteq W\). Since \(|\mathcal Q| < \kappa\),
        there is some \(A\in \mathcal Q\) such that \(H_A\subseteq W\).

        Let \(({\preceq},{\simeq})\) be such that \(A = A({\preceq},{\simeq})\).
        In other words, if \(j,k\in \sigma\) and \(\xi\in A\), then 
        \(|D_j(\xi)|_{\Bbbk} \leq |D_k(\xi)|_{\Bbbk}\) if and only if
        \(j \preceq k\) and \(D_j(\xi) = D_k(\xi)\) if and only if \(j\simeq k\).
    
        Let \(\beta = \rank({\preceq})\), so \(\beta < \theta(\sigma) < \kappa\).
        For \(\nu < \beta\), let \[\sigma_\nu = \{j\in \sigma : \rank_{\preceq}(j) = \nu\}\]
        (Note that \(\sigma_\nu\) may contain \(\simeq\)-inequivalent \(j,k\in \sigma\)
        since we are not assuming that the Ketonen order is linear.
        It will be important, however, that by the semi-linearity of the Ketonen order,
        and in particular \cref{prp:rank_closed}, 
        \(\sigma_\nu\) does not
        contain too many \(\simeq\)-equivalence classes.)
        Now
        \[\bigcap_{j\in \sigma} (F\restriction A)\ulim_{\xi < \delta} D_\xi(j) = 
        \bigcap_{\nu < \beta}\bigcap_{j\in \sigma_\nu} (F\restriction A)\ulim_{\xi < \delta} D_\xi(j)\]
        Applying \cref{clm}, 
        there is some \(\nu < \beta\) such that
        \[\bigcap_{j\in \sigma_\nu} (F\restriction A)\ulim_{\xi < \delta} D_\xi(j)\subseteq W\]

        For \(j\in \sigma_\nu\), let 
        \(H_j = (F\restriction A)\ulim_{\xi < \delta} D_\xi(j)\), so that
        \(\bigcap_{j\in \sigma_\nu} H_j\subseteq W\).
        Let \(I = \sigma_\nu/{\simeq}\) be the set of equivalence classes
        of \(\sigma_\nu\) modulo \(\simeq\).
        If \(j,k\in \sigma_\nu\) and \(j \simeq k\), then
        \(D_\xi(j) = D_\xi(k)\) for all \(\xi\in A\),
        and so \(H_j = H_k\).
        For each \(\simeq\)-equivalence class \(x\in I\),
        one can therefore define \(H_x\) to be the common value of \(H_j\) for all \(j\in x\).

        Then
        \[\bigcap_{x\in I} H_x =  \bigcap_{x\in I} \bigcap_{j\in x} H_x
                                = \bigcap_{j\in \sigma_\nu} H_j \subseteq W \]
        On the other hand, if \(x\in I\), then for any \(j\in x\), we have 
        \[F_j = F\ulim_{\xi < \delta} D_\xi(j) \subseteq 
        (F\restriction A)\ulim_{\xi < \delta} D_\xi(j) = H_j = H_x\]
        and so since \(F_j\) does not extend to a \(\kappa\)-complete ultrafilter, 
        neither does \(H_x\).
        In particular, \(H_x\nsubseteq W\).

        Fix \(\xi\in A\), and note that the set \(I\)
        is in bijection with the set
        \(\mathcal D = \{D_\xi(j) : j\in \sigma_\nu\}\). Since
        any two embeddings in \(\sigma_\nu\) have the same rank in \(\preceq\), any
        two ultrafilters in \(\mathcal D\) have the same
        rank in the Ketonen order. 
        Applying \cref{prp:rank_closed} (and the comments following \cref{def:X_closed_rank_Berkeley}), 
        the least rank Berkeley cardinal
        is \(\mathcal D\)-closed rank Berkeley, and hence it is \(I\)-closed rank Berkeley
        since \(|I| = |\mathcal D|\). 
        
        By our choice of \(\kappa\), it follows that \(\kappa\) is \(I\)-closed
        almost extendible.
        It follows from the argument of \cref{prp:normal_ultrafilters} 
        that there is an \(I\)-closed fine filter \(\mathcal Z\) on \(P(P(\delta))\)
        concentrating on the set of \(\tau\subseteq P(\delta)\) such that \(\dscott(\tau) < \kappa\).
        Since \(\mathcal Z\) is fine, for each \(x\in I\), for \(\mathcal Z\)-almost all \(\tau\),
        there is some \(A\in \tau\) such that \(A\in H_x\setminus W\).
        (This is just because for each \(\tau\), there is some \(A\subseteq \delta\) such that
        \(A\in H_x\setminus W\), and fineness implies that the set 
        \(\{\sigma\subseteq P(\delta) : A\in \sigma\}\in \mathcal Z\).)
        Since \(\mathcal Z\) is \(I\)-closed,
        these quantifiers can be exchanged: 
        for \(\mathcal Z\)-almost all \(\tau\),
        for all \(x\in I\), there is some \(A\in \tau\) such that \(A\in H_x\setminus W\).
        Therefore fix such a \(\tau\) with \(\dscott(\tau) < \kappa\).
        Then \(\bigcup (\tau\setminus W)\) belongs to \(\bigcap_{x\in I}H_x\),
        but it does not belong to \(W\) since \(W\) is \(\kappa\)-complete. 
        (Here we use \cref{lma:ordinal_complete} again.)
        This contradicts the fact that \(\bigcap_{x\in I} H_x\subseteq W\).
        This contradiction establishes that \(G\cup F\) generates a filter, completing
        the proof of the claim.
    \end{proof}
\end{thm}

\section{Questions}
For any ordinal \(\delta\), let \(F_{\delta}\) be the filter on \([\delta]^\omega\) 
generated by sets of the form \([C]^\omega\) where \(C\) is
\(\omega\)-closed unbounded in \(\delta\). 
\begin{qst}
    Assume there is a rank Berkeley cardinal.
    Is there a regular cardinal \(\delta\) such that
    \(F_{\delta}\) is atomic?
\end{qst}
Given the atomicity of the \(\omega\)-closed unbounded filter,
one would expect to prove such an analog of the partition property \(\delta\to (\delta)^\omega\), 
but the techniques of this paper seem to be powerless in the
face of a filter on a set that cannot be wellordered.

\begin{qst}
    Assume \(\lambda\) is a rank Berkeley cardinal.
    Must there be a weakly Mahlo cardinal above \(\lambda\)?
\end{qst}
\bibliographystyle{plain}
\bibliography{Bibliography}
\end{document}